\documentclass[11pt,leqno]{article}

\usepackage{amssymb,amsmath,amsthm}

\theoremstyle{plain}
\newtheorem{thm}{Theorem}
\newtheorem{pro}{Proposition}
\newtheorem*{slem}{Sublemma}

\theoremstyle{definition}
\newtheorem*{clm}{Claim}

\theoremstyle{remark}
\newtheorem*{rk}{Remark}

\newcommand{\n}{\noindent}
\newcommand{\vp}{\varepsilon}
\newcommand{\cl}[1]{\mathcal{#1}}
\newcommand{\bb}[1]{\mathbb{#1}}

\theoremstyle{plain}
\newtheorem{lem}[thm]{Lemma}

\newtheorem{pp}[thm]{Proposition~A\!\!}
\def\CC{\bb C}

\def\ovl{\overline}

\def\CC{\bb C}

\def\NN{\bb N}

\begin{document}

\title{Completely bounded maps into certain Hilbertian operator spaces}

\author{by\\ Gilles Pisier\footnote{Partially supported by   NSF      and
   Texas Advanced Research
 Program 010366-163}\\
Texas A\&M University\\
College Station, TX 77843, U. S. A.\\
and\\
Universit\'e Paris VI\\
Equipe d'Analyse, Case 186, 75252\\
Paris Cedex 05, France}

\date{}
\maketitle

\abstract{
We prove a factorization of completely bounded maps from a $C^*$-algebra
$A$ (or an
exact operator space $E\subset A$) to $\ell_2$ equipped with the operator
space structure
of   $(C,R)_\theta$ ($0<\theta<1$) obtained by complex
interpolation between the column and row Hilbert spaces.
More precisely, if $F$ denotes
$\ell_2$ equipped with the operator space structure
of   $(C,R)_\theta$, then $u:\ A \to F$ is completely bounded
iff there are states $f,g$ on $A$ and $C>0$  such that
\[ \forall a\in A\quad \|ua\|^2\le C f(a^*a)^{1-\theta}g(aa^*)^{\theta}.\]
This extends the case
$\theta=1/2$ treated in a recent paper with Shlyakhtenko \cite{PS}.
The
constants we obtain tend to 1 when $\theta \to 0$ or $\theta\to 1$, so that
we recover,
when $\theta=0$ (or $\theta=1$),
the case of mappings into $C$ (or into $R$), due to Effros and Ruan.
We use  analogues of ``free Gaussian" families in
non semifinite von Neumann algebras.
  As an application, we obtain that, if $0<\theta<1$,
 $(C,R)_\theta$ does not embed completely isomorphically
into the predual of a semifinite von Neumann algebra.
Moreover, we characterize the subspaces  $S\subset R\oplus C$ such that
the dual operator space $S^*$ embeds (completely isomorphically)
into $M_*$ for some semifinite von neumann algebra $M$: the only possibilities
are $S=R$, $S=C$, $S=R\cap C$ and direct sums built out of these three spaces.
We also discuss when $S\subset R\oplus C$ is injective, and give
a simpler proof of a result due to Oikhberg on this question.
In the appendix, we present a proof of Junge's theorem that $OH$ embeds
completely isomorphically into a non-commutative $L_1$-space. The main
idea is
similar to Junge's, but we base the argument on complex interpolation and
Shlyakhtenko's generalized circular systems
(or ``generalized free Gaussian"), that somewhat unifies
Junge's ideas with those of our work with Shlyakhtenko \cite{PS}.}
\vfill\eject

\indent\section*{Introduction}
\indent
An operator space is a Banach space given together with an isometric embedding 
into the space $B(H)$ of all bounded operators on a Hilbert space $H$.
Like the previous papers \cite{PS} and \cite{P0} (to which this one is a
natural
sequel) this paper mainly studies questions about Hilbertian operator
spaces. As
is well known, the Hilbert space $\ell_2$ can be equipped with many different
``operator space structures'', i.e.\ there are many inequivalent ways to embed
$\ell_2$ into the space $B(H)$ of all bounded operators on a Hilbert space
$H$.
Here the ``inequivalence'' is with respect to the operator space theory where
the relevant notion of morphism and isomorphism are those of
  ``completely bounded map'' and ``complete
isomorphism''.  The basic theory of operator spaces 
and their  duality was developed by Effros--Ruan and 
Blecher--Paulsen (cf.\ \cite{ER2,P3}).

The two
basic ways to realize $\ell_2$ as an operator space are the ``row'' and
``column'' ways that are defined respectively as follows. Let
\begin{align*}
R &= \overline{\text{span}}\{e_{1j}\mid j\ge 1\}\subset B(\ell_2)\\
C &= \overline{\text{span}}\{e_{i1}\mid i\ge 1\}\subset B(\ell_2).
\end{align*}
Then $R$ and $C$ are isometric to $\ell_2$ but are not completely isomorphic
(see e.g.\ \cite[p. 21]{P3}).

 In 
\cite{P1}, the author introduced another isometric embedding of $\ell_2$ into 
$B(H)$ and denoted by $OH$ the space $\ell_2$ equipped with this operator space 
structure. Being self-dual, the space $OH$ appears as the natural analogue of $\ell_2$ in 
operator space theory. It also appears as a ``midpoint'' between $R$ and $C$ in 
the sense of interpolation theory. More precisely, we have $OH = 
(C,R)_{\frac12}$. The latter formula roughly means that $OH$ is obtained by some 
sort of ``geometric mean'' out of the pair $(C,R)$. 

In \cite{P2}, the author introduced more generally a ``natural'' operator space 
structure on any $L_p$-space $(1\le p\le \infty)$. When $p=2$, one recovers the 
space $OH$ (up to cardinality considerations). This uses
 the theory   developed by Effros--Ruan and 
Blecher--Paulsen. Roughly, one starts with the obvious  (by 
Gelfand theory) operator space structure on von Neumann or $C^*$-algebras, then 
by operator space duality one passes to their duals and preduals and finally by 
interpolation one obtains the ``natural'' operator space structure on 
non-commutative $L_p$-spaces. In Banach space theory, the fact that $\ell_2$ 
embeds in $L_p = L_p([0,1], dt)$ for any $0<p<\infty$ plays a very important 
role. The case $p=1$ is probably the most useful one;  it is 
closely related to Grothendieck's theorem. More precisely, define $J\colon\ 
\ell_2\to L_p$ by setting $Jx = \Sigma x_ng_n$ where $(g_n)$ is an independent 
identically distributed (i.i.d.\ in short) sequence of  
standard (complex valued) Gaussian random 
variables.
 If we normalize $g_n$ so that $\|g_1\|_p=1$ (note $\|g_n\|_p = 
\|g_1\|_p$), then $J$ is an isometric embedding. If we replace $(g_n)$ by the 
Rademacher functions $(r_n)$, then by Khintchine's inequalities, $J$ becomes an 
isomorphic embedding for any $0<p<\infty$.

In operator space theory, the question whether $OH$ analogously embeds 
completely isomorphically (or even hopefully completely isometrically) 
into a non-commutative $L_p$-space arose immediately after \cite{P2}. Curiously, the above 
embedding $J$ produces a different operator space structure than that of $OH$. 
The latter structure actually depends on $p$ (with a sharp distinction between 
the cases $p<2$ and $p>2$) and gives $OH$ only when $p=2$. Moreover, if we 
replace $(g_n)$ or $(r_n)$ by their free analogue in Voiculescu's sense (see 
\cite{VDN}), the resulting operator space remains the same (completely 
isomorphically) for each $p$  as long as $p<\infty$. These facts are closely related to 
Lust--Piquard's non-commutative Khintchine inequalities (see \cite{P2} 
pp.~104--112 and pp.~115--121 for full details).
Passing to the Fermionic analogue of Gaussian variables (see \cite{J0}) also 
produces the same result; thus neither  does this yield a realization of $OH$ in 
non-commutative $L_p$. Actually, Junge had observed early on that
(by the Lust-Piquard inequalities)  
such a realization  is 
impossible for $2<p<\infty$ (\cite{J0} and see also the final remark below), but the case $1\le 
p<2$ and especially the crucial case $p=1$ resisted all efforts until recently 
when Junge proved that $OH$ does embed completely isomorphically into $M_*$ for some 
non-semifinite von Neumann algebra $M$. Shortly before Junge's paper, 
Shlyakhtenko and the author \cite{PS} already had used free analogues of Gaussian 
variables but in the
  {\em non-semifinite\/}  case to prove an operator space 
version of Grothendieck's theorem. With these two papers \cite{PS,J} the 
{\em non-tracial\/}  or {\em non-semifinite\/}
  theory  (also referred to as ``type III'') made its ``d\'ebut'' on 
the operator space scene. Subsequently, the author \cite{P0} showed that no 
embedding $OH\subset M_*$ can exist if $M$ is semifinite.

 This paper continues 
this line of research. We first characterize completely bounded maps from a 
$C^*$-algebra (or an exact operator space) into $OH$ (or into the space that we 
denote by $(C,R)_\theta$), thus refining \cite{PS}.

Junge actually announced (see \cite{JX1}) that his method yields
that all subspaces of quotients of
the direct sum $R\oplus C$ (of which  $OH$ is an example) 
also embed into $M_*$ for some $M$. This naturally led
us to investigate which  (infinite dimensional) quotients of
  $R\oplus C$ can embed into $M_*$ if $M$ is semifinite.
  The answer 
(see Theorem 6 below) is quite satisfactory, there are mainly 3 cases:\ $R,C$ (these embed 
into the trace class) and the space (that is denoted by $R+C$ in \cite{P1} and 
\cite{P2}) corresponding to the range of the above embedding $J$ when $p=1$. There are 
4 other possibilities built out of the first 3 using direct sums and we prove that these 7 
possibilities exhaust the list. Finally, in the appendix we show that Junge's 
embedding for $OH$ (or for subspaces of quotients of
  $R\oplus C$)  can be obtained rather quickly using Shlyakhtenko's generalized free 
Gaussian random variables. This new proof, based
on complex interpolation,   is  much quicker  and yields
  a better constant,  
  but the completely isometric case remains open.
 After this proof circulated Junge and Xu found 
another   proof using real interpolation (see \cite{JX1}).

We refer either to \cite{ER2},
\cite{Pa} or to \cite{P3} for background on operator spaces and completely
bounded maps.
 As usual, we will abbreviate ``completely bounded" by c.b. and we will
say that two operator spaces $E,F$ are completely $c$-isomorphic if there
is an isomorphism $u\colon\ E\to F$ with $\|u\|_{cb} \|u^{-1}\|_{cb}\le
c$. If $E$ is a subspace of $F$, we will say that
$E$ is completely $C$-complemented in $F$ if there is a projection 
$P\colon\  F \to E$ with $\|P\|_{cb}\le C$.
  We  recall that, given Hilbert spaces $H, K$,
 the ``minimal" (or ``spatial")
tensor product of two operator spaces $E\subset B(H)$ and
 $F\subset B(K)$ is denoted by $E\otimes_{\min} F$, it is naturally embedded
 in $B(H \otimes_2 K)$ and its norm is denoted by $\|.\|_{\min}.$
The Hilbert-Schmidt norm of a mapping $u\ : H\to K$ will be denoted by $\|u\|_2$.
We will use several times the well known fact (cf. e.g. \cite[p. 21]{P3}) that for any
mappings $u\ : R\to C$ and   $v\ : C\to R$, we have
\begin{equation}\label{RC2}
\|u\ : R\to C\|_{cb}=\|u\|_2 \ \ {\rm and}\ \  \|v\ : C\to R\|_{cb}=\|v\|_2,
\end{equation}
while for mappings from $C$ to itself or from $R$ to itself, the cb-norm
coincides with the operator norm.

\n More generally, for any Hilbert space $H$, similar properties hold for
 the associated 
``column" (resp.\ ``row") operator spaces
denoted by $H_c$ (resp.\ $H_r$).  These are 
  defined by $H_c = B({\bb C},H)$ (resp.\ $H_r = B(H^*,{\bb C})$).
 When $H=\ell_2$, 
then $H_c=C$ and $H_r=R$. $H_c$ and $H_r$ are nothing 
but analogs of $C$ and $R$ relative to general cardinals instead of that of ${\bb N}$.
When $H$ is $n$-dimensional, we denote $H_c$ (resp.\ $H_r$) by 
$C_n$ (resp.\ $R_n$).

\n The letters WEP stand for Lance's``weak expectation property": a $C^*$-algebra $A$ has the
WEP  if the inclusion $A\to A^{**}$ factors completely contractively through
$B(H)$ (see e.g. \cite{P3} for examples of its use in operator space theory).
\medskip

\indent\section*{Main Results}
\indent
Let $E\subset A$ be an operator space, given as a closed subspace of a
$C^*$-algebra $A$. Let $u\colon \ E\to \ell_2$ be a linear map. We will
identify
$\ell_2$ successively with $R,C$, and other operator spaces isometric to
$\ell_2$.

By \cite{ER1}, we know that $\|u\colon \ E\to R\|_{cb}\le 1$ 
 iff there is a state $f$   on
$A$
such that   \begin{equation}\label{eq011} 
 {\forall x\in E}\quad  \|ux\|^2 \le f(xx^*). \end{equation}
Equivalently, by \cite{ER1} this holds iff for any finite sequence
 $x_1,\ldots, x_n$ in $E$ we have
\begin{equation}\label{eq01}\sum\|ux_i\|^2 \le \left\|\sum x_ix^*_i\right\|.
\end{equation}
Similarly, $\|u\colon \ E\to C\|_{cb}\le 1$  iff there is a state $g$   on
$A$
such that   \begin{equation}\label{eq022}  
{\forall x\in E}\quad  \|ux\|^2 \le f(x^* x).   \end{equation}
Moreover, this holds iff 
for any finite sequence 
$x_1,\ldots, x_n$ in $E$ we have
\begin{equation}\label{eq02}\sum\|ux_i\|^2 \le \left\|\sum x^*_ix_i\right\|.
\end{equation}

In \cite{P1}, the author introduced a different operator space
structure on $\ell_2$, namely the space $OH$, an operator space isometric to
$\ell_2$ and uniquely characterized among operator spaces by the property that
it
is (canonically) completely isometric to its anti-dual.

For any $0<\theta<1$, one can extend (see \cite{P1}) complex interpolation to
the operator space context. Applied to the interpolation pair $(C,R)$
(using the
transposition map $x\to {}^tx$ to define ``compatibility'' in the
interpolation
sense) this method produces new operator spaces, denoted by $(C,R)_\theta$,
that
are each isometric to $\ell_2$. For $\theta=1/2$ we recover the space $OH$. To
abbreviate we will denote simply $R[\theta] = (C,R)_\theta$.
With this notation we have $R[1/2] = (C,R)_{1/2} = (R,C)_{1/2} = OH$,
and also  $R[\theta]^* =R[1-\theta]$ completely
isometrically. By convention, we set $R[0]=C$ and $R[1]=R$.

The operator space structure on $R[\theta]$ can be described more
explicitly as
follows. Let us denote by $\{e_i(\theta)\mid i=1,2,\ldots\}$ an orthonormal
basis in $R[\theta]$ (recall $R[\theta]\simeq \ell_2$ as Banach space).
Then for
any finite sequence $(a_i)$ in $B(\ell_2)$ we have
\[
\|\Sigma a_i\otimes e_i(\theta)\|_{B(\ell_2)\otimes_{\text{min}} R[\theta]} =
\sup\{(\Sigma\|s^\theta a_it^{1-\theta}\|^2_2)^{1/2}\mid s\ge 0, t\ge 0,
\|s\|_2
\le 1, \|t\|_2 \le 1\}
\]
where $\|~~~\|_2$ denotes the Hilbert--Schmidt norm. 
Equivalently,
let $p=(1-\theta)^{-1}$ and $p'=\theta^{-1}$; let $S_p$ and $S_{p'}$
denote the corresponding Schatten classes. Then the left side is equal to
the norm of the mapping $x\to \sum a_i^* x a_i$ on $S_{p'}$ and also equal to
that of the mapping $y\to \sum a_i y a_i^*$ on $S_{p}$.

In the extreme case $\theta=0$ (resp.\ $\theta=1$) we recover the space $C$
(resp.\ $R$) and the above supremum is equal to $\|\Sigma a^*_ia_i\|^{1/2}$
(resp.\ $\|\Sigma a_ia^*\|^{1/2}$). The space $R[\theta]$ 
can also be
described
as the space of ``row matrices'' inside the Schatten class $S_{p}$ with $p =
1/(1-\theta)$, when the latter is equipped with its ``natural'' operator space
structure defined (by interpolation) in \cite{P1}. Similarly
$R[\theta]$ can   be
described
as the space of ``column matrices'' inside the Schatten class $S_{p'}$. 

 In \cite{PS} it is proved
that, if $E$ is exact, then $\|u\colon \ E\to OH\|_{cb}<\infty$ iff there is a
constant $C$ and states $f,g$ such that for all $x$ in $E$ we have
\[
\|ux\|^2 \le Cg(x^*x)^{1/2} f(xx^*)^{1/2}.
\]
The first goal of this note is to prove this result with $R[\theta]$ instead of
$OH$.
Although the ingredients are the same as in \cite{PS}, our proof is somewhat
more direct. Moreover, we are able to recover the extreme cases $\theta=0$ and
$\theta=1$ described  above (due to Effros and Ruan \cite{ER1}). Note that no
assumption on $E$ is needed in the latter extreme cases, but some assumption
(such as exactness) is definitely needed when $0<\theta<1$ (see the  remark
p.~210 in \cite{PS}).

In an appendix to this note, we present a simpler proof of Junge's recent
remarkable embedding theorem of $OH$ (or $R[\theta]$ for $0<\theta<1$) into
the
predual of a von~Neumann algebra $M$. Combined with the results of the present
note, the argument of \cite{P0} shows that such an embedding is impossible,
for
any $0<\theta<1$, if $M$ is semifinite.

Let $0<\theta < 1$.   Following   \cite[\S 2.7]{P1}, we can
view
$R[\theta]$ as  an operator space such that we have (isometrically)
$$M_n(R[\theta]) = (M_n(C), M_n(R))_\theta$$
for any $n\ge 1$. We set
$$c(\theta) = (\theta^\theta(1-\theta)^{1-\theta})^{-1}.$$

\begin{thm}\label{thm1}
Let $A$ be a $C^*$-algebra. Then for any complete contraction $u\colon \ A\to
R[\theta]$, there are states $f,g$ on $A$ such that
\begin{equation}\label{eq1}
\forall\ a\in E\qquad \|ua\| \le c(\theta) f(a^*a)^{\frac{1-\theta}2}
g(aa^*)^{\frac{\theta}2}.
\end{equation}
\end{thm}

We will use the following known fact (see \cite{P1} and \cite{P4} for
related results)
\begin{lem}\label{lem1} For any $C^*$-algebra  $B$ with the WEP and
$0<\theta<1$, we have
\[
B\otimes_{\rm min} R[\theta] = (B\otimes_{\rm min} C, B\otimes_{\rm min}
R)_\theta.
\]
with equal norms.
\end{lem}

\begin{proof}
Indeed, this is equivalent to the validity of the following isometric
identities
for any $n\ge 1$

\begin{equation}\label{eq010}
B\otimes_{\rm min} (C_n,R_n)_\theta = (B\otimes_{\rm min} C_n, B \otimes_{\rm
min} R_n)_\theta.
\end{equation}
To verify \eqref{eq010} we first observe that the case $B=K(H)$ follows from
the
definition of interpolation. Then taking the bidual of both sides of
\eqref{eq010} (still restricted to the case $B=K(H)$),  we obtain the case
$B=B(H)$. Finally, if $B$ is WEP the inclusion
$B\to B^{**}$ factors completely contractively through $B(H)$, so that (using
\eqref{eq010} for $B=B(H)$) we have a complete contraction
\[
B\otimes_{\rm min}(C_n, R_n)_\theta \to (B^{**} \otimes_{\rm min} C_n, B^{**}
\otimes_{\rm min} R_n)_\theta.
\]
But $B^{**} \otimes_{\rm min} C_n = (B\otimes_{\rm min} C_n)^{**}$ and $B^{**}
\otimes_{\rm min} R_n = (B\otimes_{\rm min} R_n)^{**}$, and the norm
induced on
$B \otimes (C_n, R_n)_\theta$ by the space
\[
((B\otimes_{\rm min} C_n)^{**}, (B\otimes_{\rm min} R_n)^{**})_\theta
\]
coincides with the norm of $(B\otimes_{\rm min} C_n, B\otimes_{\rm min}
R_n)_\theta$. Indeed (see e.g.\ \cite{P3} p.  57 for more details),
since the sets appearing on both sides are obviously identical
(and each identical to $(B^*)^{n}$), we have
  isometrically  
$$(B\otimes_{\rm min} C_n, B\otimes_{\rm min}
R_n)_\theta^*=((B\otimes_{\rm min} C_n)^*, (B\otimes_{\rm min}
R_n)^*)_\theta,$$ and hence repeating the same argument for the duals, we
have
 isometrically
$(B\otimes_{\rm min} C_n, B\otimes_{\rm min}
R_n)_\theta^{**}=((B\otimes_{\rm min} C_n)^{**}, (B\otimes_{\rm min}
R_n)^{**})_\theta.$
 Therefore, we find a
completely contractive inclusion
\[
B \otimes_{\rm min} (C_n, R_n)_\theta\to (B\otimes_{\rm min} C_n,
B\otimes_{\rm
min} R_n)_\theta.
\]
On the other hand, the fact that the converse inclusion is completely
contractive is obviously true in general (without any assumption on $B$), as
follows easily by considering an embedding of $B$ onto $B(H)$.
\end{proof}

\begin{proof}[Proof of Theorem \ref{thm1}]
We will use the following formula valid for any pair $\alpha_0,\alpha_1$ of
positive numbers:
\begin{equation}\label{eq111}
\alpha^{1-\theta}_0 \alpha^\theta_1 = \inf_{\lambda>0} \{(1-\theta)
\lambda^\theta\alpha_0 + \theta \lambda^{-(1-\theta)}\alpha_1\}.
\end{equation}
Using this, \eqref{eq1} can be rewritten as
\[
\forall~\lambda>0\qquad \forall~a\in E\qquad\qquad \|ua\|^2 \le (c(\theta))^2
\{(1-\theta)\lambda^\theta f(a^*a) + \theta\lambda^{-(1-\theta)} g(aa^*)\}.
\]
By the Hahn--Banach theorem (cf.\  e.g.\  \cite[Lemma 3.4]{H}   for
details), it suffices to show  that for all finite sequences $(a_i)$ in
$E$ and all numbers
$\lambda_i>0$, we  have
\begin{equation}\label{eq1'}
\Sigma\|ua_i\|^2 \le c(\theta)^2 \{(1-\theta)\| \Sigma \lambda^\theta_i
\  a^*_i a_i\| + \theta\|\Sigma\lambda^{-(1-\theta)}_i
a_i a^*_i\|\}.  
\end{equation}
We will use the ``generalized circular elements'' introduced in
\cite{S1},  following Voiculescu's work. Since we follow closely the
ideas in \cite{PS}, we  will be brief. Let $H$ be a Hilbert space. We
assume given a set $I$ such that
$H$ has an orthonormal basis formed of the disjoint union
\[
\{e_i\mid i\in I\} \cup \{e'_i\mid i\in I\}.
\]
Let $\cl F$ be the Full Fock space over $H$, i.e.
\[
\cl F = \bb C\oplus H \oplus H^{\otimes 2}\oplus \cdots~.
\]
 Let $\Omega$ be the unit of $\bb C$, viewed as an element in $\cl F$. For any
$h$ in $H$, we denote by $\ell(h)$ (resp.\ $r(h)$) the left (resp.\ right)
creation operator on $\cl F$ i.e.\ $x\to h\otimes x$ (resp.\ $x\to x\otimes
h$). Moreover, we set $\ell_i = \ell(e_i)$ (resp.\ $r_i = r(e_i)$) and
$\ell'_i
= \ell(e'_i)$ (resp.\ $r'_i = r(e'_i)$). We define
$$
x_i = (1-\theta) \lambda^{\theta/2}_i \ell_i + \theta
\lambda^{-(1-\theta)/2}_i
\ell'{}^*_i$$
  {{and}}
$$y_i = (1-\theta) \lambda^{(1-\theta)/2}_i r'_i + \theta \lambda^{-\theta/2}_i
r^*_i.$$
Let $\cl L$ (resp.\ $\cl R$) be the von Neumann algebra generated by
$(x_i)_{i\in I}$ (resp.\ $(y_i)_{i\in I}$) in $B(\cl F)$. Note that $\cl L$
and
$\cl R$ commute with each other.

\n Let $(a_i)$ be a finite sequence in $A$. By a well known argument (as in
\cite[p. 202]{PS}) we have
\begin{align*}
\left\|\Sigma x_i\otimes a_i\right\|_{\min} &\le (1-\theta) \left\|\Sigma
\lambda^{\theta/2}_i
\ell_i\otimes a_i\right\|_{\min} + \theta\left\|\Sigma
\lambda^{-(1-\theta)/2}_i
\ell'{}^*_i \otimes a_i\right\|_{\min}\\
&\le (1-\theta) \|\Sigma a^*_ia_i\lambda^\theta_i\|^{1/2} + \theta\|\Sigma
a_ia^*_i\lambda^{\theta-1}_1\|^{1/2}\end{align*}
 {hence by Cauchy--Schwarz}$$\le \{(1-\theta)\left\|\Sigma
a^*_ia_i
\lambda^\theta_i\right\| + \theta\left\| \Sigma a_i
a^*_i\lambda^{\theta-1}_i\right\|\}^{1/2}.$$
Therefore the proof of \eqref{eq1'} (and thus of Theorem \ref{thm1}) is
reduced
to that of the following sublemma.
\end{proof}

\begin{slem}
\begin{equation}\label{eq2}
(\Sigma\|ua_i\|^2)^{1/2} \le c(\theta) \|\Sigma x_i\otimes a_i\|_{\min}.
\end{equation}
\end{slem}

\begin{proof}
By \cite{PS} we know that $\cl L$ is QWEP, i.e.\ there is a $C^*$-algebra $B$
with the weak expectation property (WEP in short) and an ideal ${\cl
I}\subset B$
such that $\cl L\simeq B/{\cl I}$. Since $\|u\|_{cb}\le 1$, $I\otimes u$
defines a
contraction from $B\otimes_{\rm min} A$ to $B \otimes_{\rm min} R[\theta]$.
Since
$R[\theta]$ has the completely contractive approximation property, we may
clearly
assume that $u$ has finite rank. In that case, if we denote by $q\colon \
B\otimes_{\rm min} A\to (B/{\cl I}) \otimes_{\rm min} A$ the canonical map,
we must
have $(I\otimes u) \ker(q) \subset {\cl I} \otimes_{\rm min} R[\theta]$
(see e.g.\
\cite[Th.\ 15.11]{P3} for details). Therefore, $I\otimes u$ defines a
contractive  map from $(B\otimes_{\rm min} A)/\ker(q)$ to $(B\otimes_{\rm
min} R[\theta])/{\cl I}
\otimes_{\rm min} R[\theta]$. Thus
\[
\|I\otimes u\colon \ \cl L\otimes_{\rm min} A \to (B\otimes_{\rm min}
R[\theta])/({\cl I} \otimes_{\rm min} R[\theta])\| \le 1.
\]
But now since $B$ has the WEP,
by Lemma \ref{lem1} the following
isometric identity holds:
\[
B\otimes_{\rm min} R[\theta] = (B\otimes_{\rm min} C, B\otimes_{\rm min}
R)_\theta.
\]
Hence we obtain a natural contractive map
\[
(B\otimes R[\theta])/({\cl I} \otimes_{\rm min} R[\theta] )\to (\cl
L\otimes_{\rm
min} C,
\cl L\otimes_{\rm min} R)_\theta.
\]
Thus to conclude it suffices to prove the following.

\begin{clm}
Let $X_\theta = (\cl L\otimes_{\rm min} C, \cl L \otimes_{\rm min} R)_\theta$.
Then we have
\[
(\Sigma \|ua_i\|^2) \le c(\theta) \|\Sigma x_i\otimes ua_i\|_{X_\theta}.
\]
Let $(z_i)_{i\le n}$ be a finite sequence in $(C,R)_\theta = R[\theta]$. We
will
show more generally that
\begin{equation}\label{eq3}
(\Sigma\|z_i\|^2)^{1/2} \le c(\theta)\|\Sigma x_i\otimes z_i\|_{X_\theta}.
\end{equation}
\end{clm}
Let us denote by $L^2_c$ (resp.\ $L^2_r$) the completion of $\cl L$ for
the norm
$x\to \|x\Omega\|$ (resp.\ $x\to \|x^*\Omega\|$). Similarly we denote
by $R^2_c$ (resp.\ $R^2_r$) the completion of
 $\cl R$ for
the norm
$x\to \|x\Omega\|$ (resp.\ $x\to \|x^*\Omega\|$).

Note that actually $L^2_c$ (resp.\ $L^2_r$) is clearly isometric to $H$ and
the
map $x\to x\Omega$ injects $\cl L$ (resp.\ $\cl R$) into $H$.
We will denote
\[  {\cl L}_\theta=(L^2_c,L^2_r)_\theta \  \ {\rm and} \   \  {\cl
R}_\theta=(R^2_c,R^2_r)_\theta .
\]
 Clearly, for any $(b_1,\ldots,
b_n)$ in ${\cl L}$ we have
$$
(\Sigma \|b_k\|^2_{L^2_c})^{1/2}  \le \|\Sigma b^*_kb_k\|^{1/2} = \|\Sigma b_k
\otimes e_{k1}\|\\
$${and} $$
(\Sigma\|b_k\|^2_{L^2_r})^{1/2}  \le \|\Sigma b_kb^*_k\|^{1/2} = \|\Sigma b_k
\otimes e_{1k}\|
$$
hence we have a contractive inclusion
\[\quad \text{from}\quad
({\cl L}\otimes_{\rm min} C, {\cl L} \otimes_{\rm min} R)_\theta
\quad \text{to}\quad
(\ell^n_2(L^2_c), \ell^n_2(L^2_r))_\theta
\]
and the latter space can be classically identified (\cite{BL}) with
$\ell^n_2([L^2_c,
L^2_r]_\theta)=\ell^n_2({\cl L}_\theta)$.

Let $(e_k)$ be the canonical
basis of
$R[\theta]$ (corresponding  to $(e_{k1})$ or $(e_{1k})$). Let $z_i = \Sigma
z_i(k)e_k$, and let
\[
x(k) = \Sigma_i z_i(k)x_i.
\]
By the preceding discussion, to prove our claim \eqref{eq3} it suffices to show
\begin{equation}\label{eq4}
\Sigma\|z_i\|^2 \le c(\theta)^2 \Sigma\|x(k)\|^2_{ {\cl L}_\theta}.
\end{equation}

Actually we will show that there is equality in the above  \eqref{eq4}.
To verify this, we now introduce
\[
y(k) = \Sigma_i \overline{z_i(k)}\ y_i
\]
so that
\begin{equation}\label{eq44}
\Sigma \langle y(k)x(k) \Omega,\Omega\rangle = \Sigma_k \Sigma_i |z_i(k)|^2
\theta(1-\theta) = \theta(1-\theta) \Sigma\|z_i\|^2.
\end{equation}

We will use the fact that for any $(X_k)$ in $\cl L$
 and any $(Y_k)$ in $\cl R$
we have
\begin{equation}\label{eq5}
|\Sigma\langle Y_kX_k\Omega,\Omega\rangle| \le (\Sigma\|X_k\|^2_{{\cl
L}_\theta})^{1/2}  (\Sigma\|Y_k\|^2_{ {\cl R}_{1-\theta}} )^{1/2}.
\end{equation}
Indeed, we have
\begin{align*}
|\Sigma\langle Y_kX_k\Omega,\Omega\rangle| &\le (\Sigma\|X_k\Omega\|^2)^{1/2}
(\Sigma\|Y_k\Omega\|^2)^{1/2}\\
&= (\Sigma\|X_k\|^2_{L^2_c})^{1/2} (\Sigma\|Y_k\|^2_{R^2_r})^{1/2}\end{align*}
 {and also since $Y_kX_k = X_kY_k$}
$$|\Sigma\langle Y_kX_k\Omega,\Omega\rangle| \le
(\Sigma\|X_k\|^2_{L^2_r})^{1/2}
(\Sigma\|Y_k\|^2_{R^2_c})^{1/2},$$
hence by the bilinear interpolation theorem (\cite{BL}) \eqref{eq5} follows
(since $(R^2_r,R^2_c)_{\theta} = (R^2_c,R^2_r)_{1-\theta}={\cl
R}_{1-\theta}$). 

We will show that
\begin{equation}\label{eq6}
(\Sigma\|y(k)\|^2_{{\cl
R}_{1-\theta}})^{1/2} \le (1-\theta)^\theta \theta^{1-\theta}
(\Sigma\|z_i\|^2)^{1/2}.
\end{equation}
Let $f_i$ be the ${\cl R}$-valued analytic function defined on $\bb C$ by
\[
f_i(z) = ((1-\theta)^{1-z}\theta^z)^{-1} ((1-\theta) \lambda^{z/2}_ir'_i
+
\theta \lambda^{-(1-z)/2}_i r^*_i).
\]
Note that $f_i(1-\theta) = ((1-\theta)^\theta \theta^{1-\theta})^{-1}y_i$.
Moreover, for any $(\alpha_i)$ in $\ell_2(I)$, we have
\[
\left\{\begin{array}{lll}
\|\Sigma \alpha_if_i(z)\|_{R^2_c} = (\Sigma|\alpha_i|^2)^{1/2}&\text{if}&
\text{Re}(z) = 0\\
\|\Sigma \alpha_if_i(z)\|_{R^2_r} = (\Sigma|\alpha_i|^2)^{1/2}&\text{if}&
\text{Re}(z) = 1.
\end{array}\right.
\]
Hence
\begin{equation}\label{eq7}
\|\Sigma \alpha_if_i(1-\theta)\|_{{\cl
R}_{1-\theta}}  \le (\Sigma|\alpha_i|^2)^{1/2}.
\end{equation}
But now
\begin{equation}\label{eq8}
\Sigma_i \overline{z_i(k)} f_i(1-\theta) = ((1-\theta)^\theta
\theta^{1-\theta})^{-1} y(k)
\end{equation}
hence \eqref{eq6} follows from \eqref{eq7} and \eqref{eq8}. Now combining
\eqref{eq6} with \eqref{eq5} and \eqref{eq44} we find
\begin{align*}
(1-\theta)\theta \Sigma\|z_i\|^2 &= \Sigma\langle y(k)
x(k)\Omega,\Omega\rangle\\
&\le (\Sigma\|x(k)\|^2_{{\cl
L}_{\theta}})^{1/2} (\Sigma\|y(k)\|^2_{{\cl
R}_{1-\theta}})^{1/2}\\
&\le (\Sigma \|x(k)\|^2_{{\cl
L}_{\theta}})^{1/2} {(1-\theta)}^{\theta} \theta^{1-\theta}
(\Sigma\|z_i\|^2)^{1/2}
\end{align*}
which, after a suitable division, yields \eqref{eq4}. This completes the proof
of \eqref{eq3}, of the above claim, and of the sublemma.
Note that   an obvious modification of the proof of
\eqref{eq6} shows that the converse of \eqref{eq4} also holds so
\eqref{eq4} is indeed an equality.\end{proof}

\begin{rk}
Arguing as in \cite{PS}, it is easy to check that Theorem \ref{thm1} remains
valid when $A$ is replaced by an exact subspace $E\subset A$ with exactness
constant $\le c$, provided the constant $c(\theta)$ in \eqref{eq1} is
replaced
by $c(\theta)c$.
\end{rk}

\begin{rk}
Note that $c(\theta)\to 1$ when either $\theta\to 0$ or $\theta\to 1$. In the
cases $\theta=0$ and $\theta=1$, Theorem \ref{thm1} is well known (cf.\
\cite{ER1}). In that case, Theorem \ref{thm1} still holds  when $A$ is
replaced by  an {\em arbitrary\/} subspace $E\subset A$. However, when
$0<\theta<1$, some  extra assumption (such as exactness) is necessary.
Indeed, if we take
$\theta=1/2$, let $(a_i)$ be the orthonormal basis of $OH = (C,R)_{1/2}$, and
let $u$ be the identity map, we have $\left\|\sum\limits^n_1
a^*_ia_i\right\|^{1/2} = \|\sum\nolimits^n_1 a_ia^*_i\|^{1/2} = n^{1/4}$ but
$\sum\limits^n_1 \|ua_i\|^2 = \sum\nolimits^n_1 \|a_i\|^2 = n$,
which shows that
\eqref{eq1'} fails. Similarly, the extension property valid when either
$\theta= 0$ or $\theta=1$ is no longer true in general, indeed this
is closely related to the fact that $R$ or $C$ are injective operator
spaces, while $R[\theta]$ is not when $0<\theta<1$.

\end{rk}

\begin{rk}
In the preceding argument, the only delicate point is \eqref{eq4}. Note that
actually, it is easy to show that equality holds in \eqref{eq4}. We chose to
base the above proof of \eqref{eq4} solely on complex interpolation to make it
accessible to a reader unfamiliar with the Tomita--Takesaki theory.
However, if
one uses the latter theory, in the form made explicit by Shlyakhtenko in
\cite{S1}, it is very easy to explain why \eqref{eq4} should be true. We
now  review this
alternate approach.
\end{rk}

\begin{proof}[Alternate proof of \eqref{eq4}]
 Let $\varphi$ be the vacuum state, defined on $B(\cl F)$ by $\varphi(T) =
\langle T\Omega,\Omega\rangle$. Let $\xi_i  = \theta(1-\theta)^{-1}
\lambda^{-1/2}_i$. Note that
\begin{equation}\label{a}
x_i = (1-\theta)\lambda^{\theta/2}_i (\ell_i  + \xi_i\ell'{}^*_i);
\end{equation} 
therefore $\cl L$ can be viewed as generated by $\{\ell_i
+\xi_i\ell'{}^*_i\}$.
We define a one parameter group of unitary operators $u_t$ on $H$ by
setting for
any $t$ in $\bb R$
\[
\forall~j\in I\qquad u_te_j = (\xi_j)^{2it} e_j, \quad u_te'_j  = (\xi_j)^{-2it}
e'_j.
\]
We extend $u_t$ (by the so-called first quantization) to a unitary operator
$U_t$ on $\cl F$ such that $U_t\Omega=\Omega$ and $U_t = u_t\otimes\cdots
\otimes u_t$ on $H\otimes\cdots\otimes H$. For any $x$ in $\cl L$, we denote
\[
\sigma_t(x) = U_txU^{-1}_t.
\]
Note that $\sigma_t(\ell_j) = \xi^{2it}_j\ell_j$ and $\sigma_t(\ell'_j) =
\xi^{-2it}_j\ell'_j$ so that we have for all $j$:
\[
\sigma_t(x_j) = \xi^{2it}_jx_j\quad \text{and}\quad \sigma_t(x^*_j) =
\xi^{-2it}_j x^*_j.
\]\
Therefore $\sigma_t$ is a one parameter group of automorphisms of $\cl L$,
that
is nothing but the classical modular automorphism group of $\cl L$ relative to
the state $\varphi$. In particular, $(\sigma_t)$ satisfies the KMS condition:\
for any polynomials $x,y$ in the generators $\{x_j\}$ we have
\[
\varphi(\sigma_i(x)y) = \varphi(yx)
\]
where $z\to \sigma_z$ is the obvious analytic extension of $t\to\sigma_t$.
As is well known (cf.\ e.g.\ \cite{PW1,PW2}), we have in this situation
\[
\|x\|_{{\cl L}_\theta}=
\|x\|_{(L^2_c,L^2_r)_\theta} = \|\sigma_{-i\theta/2}(x)\|_{L^2_c} =
\varphi(\sigma_{-i\theta/2} (x)^* \sigma_{-i\theta/2}(x))^{1/2}.
\]
Hence, we can write since $\xi_j\lambda^{1/2}_j = \theta(1-\theta)^{-1}$
\begin{align*}
\Sigma\|x(k)\|^2_{\cl L_\theta} &= \Sigma\varphi(\sigma_{-i\theta/2} (x(k))^*
\sigma_{-i\theta/2} (x(k)))\\
&= \Sigma\|\sigma_{-i\theta/2} (x(k)) \Omega\|^2\\
&= \sum_k \left\|\sum_j z_j(k) \sigma_{-i\theta/2}(x_j)\Omega\right\|^2\\
&= \sum_k \left\|\sum_j z_j(k) \xi^\theta_j x_j\Omega\right\|^2\\
&= \sum_k \left\|\sum_j z_j(k) \xi^\theta_j (1-\theta) \lambda^{\theta/2}_j
e_j\right\|^2\\
&= \sum_j \|z_j\|^2 (c(\theta))^{-2}.
\end{align*}
Hence, we obtain the announced equality
\[
(\Sigma\|z_j\|^2)^{1/2} = c(\theta) (\Sigma\|x(k)\|^2_{\cl L_\theta})^{1/2}.
\qquad  \qed
\]
\renewcommand{\qed}{}\end{proof}

 The converse of Theorem \ref{thm1} also holds,
as  follows. In the case $\theta=1/2$, this was proved in \cite{PS}. We
give a more direct argument.

\addtocounter{pro}{2}
\begin{pro}\label{pro2}
\addtocounter{thm}{1}
Let $E\subset A$ be an operator space embedded in a $C^*$-algebra $A$.
 Then any linear map $u\colon \ E\to \ell_2$    for which there
are states $f,g$ on $A$ and a constant $C$ such that
\begin{equation}\label{eq9}
\forall~a\in E\qquad \|ua\| \le Cf(a^*a)^{(1-\theta)/2} g(aa^*)^{\theta/2}
\end{equation}
is completely bounded, with $\|u\|_{cb} \le C$,  as a mapping into
$\ell_2$ equipped with the operator space structure of $(C,R)_\theta$.
\end{pro}

\begin{proof}
Let $F = (C,R)_\theta$ (recall this is isometric to $\ell_2$).
 We will show that for any $a = (a_{ij})$ in the unit ball of
$M_n(E)$ we have $\|(u(a_{ij}))\|_{M_n(F)} \le C$. Let $(T_k)$ be an
orthonormal basis of $F = (C,R)_\theta$. Let
$u_k\colon \ E\to \bb C$ be defined by $u(a) = \Sigma u_k(a)T_k$. Let us
denote
$\gamma_k = \sum\limits_{ij} u_k(a_{ij}) e_{ij}\in M_n$. Using $M_n(F) \simeq
M_n\otimes F$, the matrix $(u(a_{ij}))$ can then be rewritten as
\[
\Sigma e_{ij} \otimes u(a_{ij}) = \sum_k \gamma_k \otimes T_k.
\]
Let $p,p'$ be defined by
$1-\theta=1/p$ and $\theta=1/p'$.
By definition of $(C,R)_\theta$, we have (see the identity (8.5), p.~83 in
\cite{P1}, but note that our space $R[\theta]$  corresponds to the space
denoted by
$R(1-\theta)$ in \cite{P1}):
\begin{equation}\label{eq0100}
\|(u(a_{ij}))\|_{M_n(F)} = \sup\left\{\left(\sum_k \|s\gamma_k
t\|^2_2\right)^{1/2}\right\}
\end{equation}
where $\|~~\|_2$ denotes the Hilbert--Schmidt norm on $M_n$ and where the
supremum runs over all pairs $(s,t)$ in $(M_n)_+ \times (M_n)_+$ such that
${\rm tr}\
s^{2p'} \le 1$ and ${\rm tr}\ t^{2p} \le 1$.

Let $x_{ij}$ denote the $(ij)$-entry of $s at$ (so that $x_{ij} =
\sum\limits_{k\ell} s_{ik} a_{k\ell} t_{\ell j}$). We have
\[
\sum_k \|s\gamma_kt\|^2_2 = \sum_{ij} \|u(x_{ij})\|^2_F.
\]
We claim that \eqref{eq9} implies
\begin{equation}\label{eq10}
\sum_{ij} \|u(x_{ij})\|^2_F \le C^2.
\end{equation}
By \eqref{eq0100}, this claim implies that $\|u\|_{cb}\le C$, thus
completing the proof.
To show this claim, we may as well assume (replacing $a$ by $v_1 a v_2$ for
suitable  unitaries $v_1,v_2$ in $M_n$, and using \eqref{eq0100} again)
that $s$ and $t$ are diagonal matrices. We then  have
\[
\sum_{ij} \|u(x_{ij})\|^2 = \sum_{ij} s^2_{ii}\|u(a_{ij})\|^2 t^2_{jj}
\]
hence, since $\Sigma s^{2p'}_{ii}\le 1$ and $\Sigma t^{2p}_{jj} \le 1$, we
have using
\eqref{eq9}
\[
\sum_{ij} \|u(x_{ij})\|^2 \le C^2\left\|\Sigma f(a^*_{ij}a_{ij})^{1-\theta}
g(a_{ij}  a^*_{ij})^{\theta} e_{ij}\right\|_{B(\ell^n_p,\ell^n_p)}.
\]
Thus, the above claim follows from the next lemma.
\end{proof}

\begin{lem}\label{lem3}
Let $f,g$ be states on a $C^*$-algebra $A$. We have then for any $n\ge 1$ and
any $a$ in $M_n(A)$
\[
\left\|\sum_{ij} f(a^*_{ij}a_{ij})^{1-\theta} g(a_{ij}a^*_{ij})^{\theta}
e_{ij}\right\|_{B(\ell^n_p,\ell^n_p)} \le \|a||_{M_n(A)}.
\]
\end{lem}

\begin{proof}
To prove Lemma \ref{lem3}, we first check that, if $\|a\|_{M_n(A)}
\le 1$, then $\alpha_0(i,j) = f(a^*_{ij}a_{ij})$ and $\alpha_1(i,j) =
g(a_{ij}a^*_{ij})$ satisfy 
\begin{equation}\label{eq11}
\sup_j \Sigma_i \alpha_0(i,j) \le 1 \quad \text{and}\quad \sup_i \Sigma_j
\alpha_1(i,j) \le 1.
\end{equation} Indeed, we have for any fixed $j$
\[
\sum_i f(a^*_{ij}a_{ij}) \le \left\|\sum_i a^*_{ij}a_{ij}\right\| =
\left\|\sum_i a_{ij} \otimes e_{ij}\right\| \le \|a\|_{M_n(A)},
\]
and similarly for the other sum.
Now let $\alpha_0(i,j), \alpha_1(i,j)$ be 
any $n\times n$ matrices with nonnegative
entries satisfying \eqref{eq11}.
Then it is well known and easy to check by interpolation
 (see e.g. \cite{P5}   for more on this topic) that
\[
\left\|\sum_{ij}\alpha_0(i,j)^{1-\theta}
\alpha_1(i,j)^{\theta}e_{ij}\right\|_{B(\ell^n_p,\ell^n_p)}\le
1,
\]
or equivalently for any $s_i$, $t_j\ge 0$ with $\Sigma s_i \le 1$, $\Sigma
t_j\le 1$ we have
\[
\sum_{ij} \alpha_0(i,j)^{1-\theta} \alpha_1(i,j)^{\theta} s^{1/p'}_i
t^{1/p}_j \le 1.
\]
Indeed, by H\"older's inequality (recall $1/p=1-\theta$ and $1/p'=\theta$)
this is
\begin{align*}
&\le \left(\sum_{ij} \alpha_0(i,j)t_j\right)^{1-\theta} \left(\sum_{ij}
\alpha_1(i,j)s_i\right)^{\theta}\\
&\le \left(\sum_j t_j\right)^{1-\theta} \left(\sum_i s_i
 \right)^{\theta}\\
&\le 1.
\end{align*}

\end{proof}

\begin{rk}
When $\lambda_i=1$ for all $i$ and $\theta=1/2$, we have $x_i = (1/2)x'_i$ where
\[
x'_i = \ell_i + \ell'{}^*_i.
\]
Then $(x'_i)$ is a free circular (i.e. free analogue of complex Gaussian) family in 
Voiculescu's sense (cf.\ \cite{VDN}). It is easy to see in this case that for 
any finite sequence $(a_i)$  in $B(\ell_2)$ we have
\[
(1/2)\left\|\sum a_i \otimes x'_i\right\| \le \max\left\{\left\|\sum 
a^*_ia_i\right\|^{1/2}, \left\|\sum a_ia^*_i\right\|^{1/2}\right\} \le 
\left\|\sum a_i\otimes x'_i\right\|.
\]
Thus $\ovl{\text{span}}[x'_i]$ is completely isomorphic to the space $R\cap C$ 
studied in \cite{HP} (see also \cite[p. 209]{P3}). The notation $R\cap C$ comes 
from the fact that if we consider $\delta_i = e_{1i} \oplus e_{i1}$ in $R\oplus 
C$ then we have for $(a_i)$ as above
\[
\left\|\sum a_i\otimes \delta_i\right\| = \max\left\{\left\| \sum a_i 
a^*_i\right\|^{1/2}, \left\|\sum a^*_ia_i\right\|^{1/2}\right\};\]
so that
\[
\ovl{\text{span}}[\delta_i] = \{(x,{}^t x)\mid x\in {R}\}\subset
R\oplus C
\]
appears as the diagonal in $R\oplus C$. 
Let ${\cl L}$ be again the von Neumann algebra generated by $\{x'_i\}$.
 We claim (see \cite{HP}, see also \cite[p. 209]{P3})
that there is a normal c.b.\ projection $P\colon  \ {\cl L}\to 
\ovl{\text{span}}[x'_i]$ with $\|P\|_{cb}\le 2$. Indeed, let $Q\colon \ {\cl 
F}\mapsto {\cl F}$ (resp.\ $Q'\colon \ {\cl F}\to {\cl F}$) be the orthogonal 
projection onto $\ovl{\text{span}}[e_i\mid i\in I]$ (resp.\ 
$\ovl{\text{span}}[e'_i\mid i\in I]$) viewed as a subspace of $H$, itself 
embedded into ${\cl F}$ via tensors of degree 1. Then the map $P$ defined by
\begin{equation}\label{d}
\forall T\in {\cl L}\qquad P(T) =\ell(Q(T\Omega)) + \ell(Q'(T^*\Omega))^*
\end{equation}
is the announced projection (see more generally Lemma \ref{lem22} below). 
Therefore $(R\cap C)^*$ embeds completely isomorphically  into ${\cl 
L}_*$ and in the present special case ${\cl L}$ is semifinite (and
actually finite).

Consider now a family $\xi =(\xi_i)$ with $\xi_i>0$. Let $\delta^\xi_i = e_{1i} 
\oplus \xi_i e_{i1} \in R\oplus C$. Thus if $(a_i)$ is as before, we have
\[
\left\|\sum a_i\otimes\delta^\xi_i\right\| = \max\left\{\left\|\sum 
a_ia^*_i\right\|^{1/2}, \left\|\sum \xi^2_i a^*_ia_i\right\|^{1/2}\right\},
\]
and also
\begin{equation}\label{c}
(1/2)\left\|\sum a_i \otimes (\ell_i + \xi_i\ell'{}^*_i)\right\| \le \left\|\sum 
a_i\otimes \delta^\xi_i\right\| \le \left\|\sum a_i \otimes (\ell_i + \xi_i 
\ell'{}^*_i)\right\|.
\end{equation}
Thus we have completely 
isomorphically
\[
\ovl{\text{span}}[\ell_i + \xi_i\ell'{}^*_i] \simeq 
\ovl{\text{span}}[\delta^\xi_i].
\]
In particular, if $(\xi_i)$ is as in \eqref{a} then 
$\ovl{\text{span}}[x_i]\simeq \ovl{\text{span}}[\delta^\xi_i]$ (completely
isomorphically).
 Note that $\ovl{\text{span}}[\delta^\xi_i]$ can also be viewed 
as the graph of the unbounded operator $\Lambda\colon \ R\to C$ taking $e_{1i}$ 
to $\xi_ie_{1i}$, with $\text{Dom}(\Lambda)
 = \{x = \Sigma x_ie_{1i}\in 
R\mid \Sigma|\xi_ix_i|^2<\infty\}$. More precisely, if we denote
\[
G(\Lambda) = \{(x,\Lambda x)\mid x\in \text{Dom}(\Lambda)\},
\]
then we have $\ovl{\text{span}}\ \delta^\xi_i\simeq G(\Lambda)$ completely 
isometrically. (Note:\ In analogy with $R\cap C$, it would be natural to denote 
$G(\Lambda)$ by $R\cap \Lambda^{-1}(C)$ but we prefer not to use this notation.)

The next result (extending the case   $\xi_i=1\  \forall i$)  is easy to deduce from Shlyakhtenko's  \cite{S1}.
\begin{lem}\label{lem22}
Let ${\cl L}$ be the von Neumann algebra generated by the family $(x_i)$ defined 
 by \eqref{a}. Then the mapping $P\colon \ {\cl L} \to {\cl L}$ defined by \eqref{d} is a 
normal c.b.\ projection from ${\cl L}$ onto $\ovl{\text{span}}[x_i]$ with 
$\|P\|_{cb}\le 2$. In particular $\ovl{\text{span}}[x_i]$ embeds completely isomorphically 
into 
${\cl L}_*$.
\end{lem}

\begin{proof}
We   first claim that $T\to P(T)$ is c.b.\ on
$B({\cl F})$ with cb-norm $\le  2$. This is easy to see. Indeed,
consider
$P_1:\ B({\cl F})\to B({\cl F})$ and
$P_2:\ B({\cl F})\to B({\cl F})$
defined by
$$P_1(T)=\ell(Q(T\Omega))\quad\text{and}\quad
P_2(T)=\ell(Q'(T^*\Omega))^*$$
so that $P(T)=P_1(T)+P_2(T)$.
We will show that $\|P_1\|_{cb}\le 1$
and 
 $\|P_2\|_{cb}\le 1$.
Indeed, the ranges of $P_1$ and $P_2$ are
respectively
$  \overline{\text{span}}[\ell_i]  $   and $ \overline{\text{span}}[\ell'{}_i^*]    $.
Assuming $I=\NN$ for simplicity, we have
(see e.g. \cite[p. 176-177]{P3})
$$  \overline{\text{span}}[\ell_i]\simeq C\quad\text{and}\quad 
  \overline{\text{span}}[\ell'{}_i^*]\simeq R .$$
Note that we have obviously (recall $\varphi$ is the vacuum state)
$$\|P_1(T)\|=\|Q(T\Omega)\|\le \|T\Omega\|=\langle T^* T \Omega, \Omega
\rangle^{1/2} =\varphi(T^* T)^{1/2}$$  
$$\|P_2(T)\|=\|Q'(T^*\Omega)\|\le \|T^*\Omega\|=\langle T T^* \Omega, \Omega
\rangle^{1/2} =\varphi(T T^*)^{1/2}$$ Hence,  $\|P_1\|_{cb}\le 1$ and
$\|P_2\|_{cb}\le 1$ (and a fortiori $\|P\|_{cb}\le 1$) follow using 
\eqref{eq022} and \eqref{eq011}.
 Thus it suffices to prove that the restriction of $P$ to ${\cl L}$ is a projection 
onto $\ovl{\text{span}}[x_i]$.
We know (see \cite{S1}) that the map $T\to T\Omega$ is faithful (i.e. injective)
 on ${\cl L}$. 
Let $T$ be a polynomial in $x_i,x^*_i$ $ (i\in I)$. We can write a priori
\[
T\Omega = \sum t_ie_i + \sum t'_ie'_i + r
\]
where $r$ is a sum of tensors of degree $>1$. By \cite[Lemma 3.2]{S1}, we know 
that the (antilinear) map $S$ taking $T\Omega$ to $T^*\Omega$ takes $r$ to 
another sum $r'$ of tensors of degree $>1$ in ${\cl F}$. Moreover, since $ 
(\ell_i+\xi_i\ell'{}^*_i)\Omega=e_i $ and $(\ell_i + \xi_i\ell'{}^*_i)^* \Omega = 
\xi_ie'_i$, we have
\[
Se_i = \xi_ie'_i\quad \text{and}\quad Se'_i = \xi^{-1}_ie_i,
\]
and hence
\[
T^*\Omega = \sum \bar t_i\xi_ie'_i  + \sum \ovl{t'_i} \xi^{-1}_ie_i+r'.
\]
Therefore we have $Q(T\Omega) = \Sigma t_ie_i$,  $Q'(T^*\Omega) = \Sigma \bar 
t_i \xi_ie'_i$, and we finally obtain
\[
P(T) = \sum t_i\ell(e_i) + \sum t_i\xi_i \ell(e'_i)^* = \sum t_i(\ell_i + 
\xi\ell'{}^*_i) \in \text{span}[x_i].
\]
In particular, we find $P(\ell_i + \xi_i\ell'{}^*_i) = \ell_i + 
\xi_i\ell'{}^*_i$ and hence $P(x_i) = x_i$ for all $i$. This proves that 
$P_{|{\cl L}}$ is a projection from ${\cl L}$ onto $\ovl{\text{span}}[x_i]$.
\end{proof}

Note that, by \eqref{c}, if
\begin{equation}\label{eq(e)}
0 < \inf \xi_i \le \sup \xi_i < \infty
\end{equation}
then $\ovl{\text{span}}[x_i]$ (or equivalently $G(\Lambda)$) 
is again completely
isomorphic to 
$R\cap C$ and hence its dual embeds in $M_*$ for some semifinite $M$. We
will now show  that if either $\inf\xi_i = 0$ or $\sup\xi_i = \infty$,
then such an embedding 
$G(\Lambda)^*\subset M_*$ with $M$ semifinite exists if and only if we have for 
some $\vp>0$
\begin{equation}\label{eq(b)}
\sum_{i\colon  \xi_i<\vp} \xi^2_i + \sum_{i\colon \xi^{-1}_i<\vp} \xi^{-2}_i < 
\infty.
\end{equation}

\end{rk}
Let $M$ be a von Neumann algebra with predual $M_*$. As already mentioned at the end of  
\cite{P0}, Theorem \ref{thm1} and its converse (Proposition \ref{pro2}) admit 
the following corollary:\ for any $0<\theta<1$, the space $R[\theta]$ does not 
embed completely isomorphically into $M_*$ when $M$ is semifinite. It can be 
shown (see the appendix below) that $R[\theta]$ is completely isometric to a 
quotient of a subspace $S$ of $R\oplus C$. Thus $R[\theta]^*$ embeds in $S^*$,
 and hence to embed $R[\theta]^*$ into $M_*$ it suffices to embed
$S^*$ into  $M_*$.
 Indeed, Marius Junge announced that if $S$ 
is any subspace of $R\oplus C$ then $S^*$ embeds completely isomorphically into 
$M_*$ for some von Neumann algebra $M$. Let $S \subset R\oplus C$ be such a 
subspace. For convenience, let us assume that $S$ is not completely isomorphic 
to either $R,C$ or $R\oplus C$. Then Q.\ Xu (\cite{X1}) observed the fact
(presumably known to Junge) that $S$ can be rewritten (up to complete
isomorphism) as a direct sum 
$H_r\oplus \widetilde S\oplus K_c$ where $H_r,K_c$ are suitable Hilbert spaces 
equipped respectively with the row and column operator space structure, and 
where $\widetilde S \subset R\oplus C$ is the (closed) graph of a (closed) densely 
defined operator $\Lambda\colon \ R\to C$, injective (on its domain) and with 
dense range. As explained in the appendix, the fact that $(\widetilde S)^*$ 
embeds into $M_*$ for some suitable $M$ can be deduced from the basic properties 
of Shlyakhtenko's generalized free circular elements, already used in \cite{PS}. 
The typical $M$ is then not semifinite. The next result shows that this cannot 
be avoided.

\begin{thm}\label{thm3}
Let $S\subset R\oplus C$ be an arbitrary infinite dimensional subspace. Then 
there is a semifinite von Neumann algebra $M$ such that $S^*$ embeds completely 
isomorphically into $M_*$ iff $S$ is completely isomorphic to one of the spaces
\[
R, C, R\oplus C, R\cap C, R\oplus (R\cap C), C \oplus (R\cap C), R\oplus (R\cap 
C) \oplus C.
\]
\end{thm}

\begin{rk}
As is well known, we have $R^*\simeq C$ and $C^*\simeq R$, so that $R^*$ and 
$C^*$ embed (completely isometrically) in $K^*\simeq S_1$ (the trace class). 
Consequently $(R\oplus C)^*\simeq R\oplus C$ embeds completely isomorphically 
into $S_1\oplus S_1 \simeq S_1$. The space $R\cap  C$ is less trivial, but it 
was shown by Lust-Piquard and the author (see \cite{LPP} or 
\cite[p.~194]{P3}) that $(R\cap C)^*$ 
embeds completely isomorphically  into  the most classical $L_1$-space, namely
$L_1([0,1],dt ) $. Therefore, we 
do have $S^*\subset M_*$ with $M$ semifinite for any of the 7 spaces in the 
above list.
\end{rk}

Our task will now be to show that the latter list is complete.

\begin{rk}\label{rem33}
Consider a (closed) subspace $S\subset R\oplus C$. As explained above, we can 
write
\begin{equation}\label{eq100}
S\simeq H_r\oplus G(\Lambda) \oplus K_c
\end{equation}
with
\[
G(\Lambda) = \{(x,\Lambda x)\mid x\in \text{Dom}(\Lambda)\} \subset R\oplus C
\]
where $\text{Dom}(\Lambda) \subset R$ is a dense subspace and $\Lambda\colon \ 
\text{Dom}(\Lambda)\to C$ is a closed unbounded operator with zero kernel and 
dense range.
\end{rk}

By the polar decomposition of $\Lambda$ and the ``homogeneity'' of $R$ and $C$ 
(in the sense of \cite[p.~172]{P3}), we may assume that $\Lambda > 0$. Using the 
spectral theory of Hermitian operators, we can then decompose $\Lambda$ as 
$\Lambda = \Lambda_1+\Lambda_2$ with $0<\Lambda_1\le 1$ and $\Lambda_2 \ge 1$, 
and consequently we may decompose
\begin{equation}\label{eq1000}
G(\Lambda) \simeq G(\Lambda_1) \oplus G(\Lambda_2)
\end{equation}
where $\Lambda_1,\Lambda_2$ are unbounded self-adjoints of the same form as 
$\Lambda$ but in addition such that $\Lambda_1$ and $\Lambda^{-1}_2$ are bounded 
with norm $\le 1$. The key to the preceding theorem then lies in the next 
statement.

\begin{lem}\label{lem4}
Consider $\Lambda>0$ with $\|\Lambda\|\le 1$ and $\Lambda^{-1}$ unbounded. Let 
$E(\vp)$ be the spectral projection of $\Lambda$ relative to $(0,\vp)$, so that 
$0\ne \|\Lambda E(\vp)\| \le \vp$ for any $\vp>0$. Assume that there is a 
semifinite $M$ such that $G(\Lambda)^*$ embeds completely isomorphically 
into $M_*$. Then, for 
$\vp>0$ small enough, $\Lambda E(\vp)$ must be Hilbert--Schmidt.
\end{lem}

\begin{proof}
The basic idea is similar to the one in \cite{P0} but the details are more 
complicated.
By assumption, we have an embedding $j\colon \ G(\Lambda)^*\subset M_*$. Let $u 
= j^*\colon \ M\to G(\Lambda)$. We may assume that $\|u\|_{cb}\le 1$ and that 
there is a constant $c$ such that for any $n$ and any $a$ in $M_n(G(\Lambda))$ 
with $\|a\|<1$, there is $\tilde a$ in $M_n(M)$ with 
\begin{equation}\label{eq101}
\|\tilde a\| < c
\end{equation}
such that $(I\otimes u) (\tilde a) =a$. Note that $u$ is ``normal'', i.e.\ is 
$(\sigma(M,M_*), \sigma(G(\Lambda), G(\Lambda)^*))$ continuous. The map $u$ can 
clearly be rewritten as $ux = (vx,\Lambda vx)$ with 
\[
\|v\colon \ M\to R\|_{cb}\le 1\quad \text{and}\quad \|\Lambda v\colon \ M\to 
C\|_{cb}\le 1.
\]
Let $\tau$ be a semifinite faithful normal trace on $M$. Since $v$ and $\Lambda 
v$ are normal, arguing as in \cite{P0}, we find normal states $f,g$ on $M$ such 
that
\begin{align*}
\|vx\| &\le f(xx^*)^{1/2}\\
\|\Lambda vx\| &\le g(x^*x)^{1/2}
\end{align*}
for all $x$ in $M$.

We may view $f,g$, as elements of $L_1(\tau)$, i.e.\ positive unbounded 
operators affiliated to $M$ such that $\tau(f) = \tau(g) = 1$, and consequently 
we will write $f(\cdot) = \tau(f\cdot)$ and $g(\cdot)  = \tau(g\cdot)$.
Fix $\alpha>0$. Let $p$ (resp.\ $q$) be the spectral projection of $f$ (resp.\ 
$g$) associated to $(\alpha^{-1},\alpha]$, so that in $M_*$, we have
 $\alpha^{-1} p\le pfp\le \alpha p$ and $\alpha^{-1} q\le q g q 
\le \alpha q$. Choosing $\alpha = \alpha(\delta)$ large enough we can ensure 
that moreover
\[
\|f(1-p)\|_{M_*} < \delta^2\quad \text{and}\quad \|g(1-q)\|_{M_*} < \delta^2.
\]
Moreover, we have
$ \alpha^{-1} \tau(p)\le \tau(f)=1$   {and} $\alpha^{-1} \tau(q)\le \tau(g)=1$,
and hence
\[ \tau(p)\le \alpha   \quad \text{and}\quad  \tau(q)\le \alpha.
\]
We then define
\[
v_\delta x = v(pxq).\]
Note that 
\[
\|v_\delta x\| \le f(pxqx^*p)^{1/2} \le \sqrt\alpha\ \tau(pxqx^*)^{1/2}  = 
\sqrt\alpha \tau(qx^*px)^{1/2} \le \sqrt\alpha \tau(qx^*x)^{1/2},
\]
where the last equality holds by  the trace property, so that by   \eqref{eq022}
\begin{equation}\label{eq102}
\|v_\delta\colon \ M\to C\|_{cb} \le (\alpha\tau(q))^{1/2}\le \alpha.
\end{equation}
On the other hand, we have
\begin{equation}\label{eq1021}
vx - v_\delta x = v_1x+v_2x
\end{equation}
with $v_1x = v(x(1-q))$ and $v_2x  = v((1-p)xq)$. Note that
\[
\|\Lambda v_1x\|\ \le g((1-q) x^*x(1-q))^{1/2} = \tau(g(1-q) x^*x)^{1/2}
\]
hence by   \eqref{eq022}  
\begin{equation}\label{eq103}
\|\Lambda v_1\colon \ M\to C\|_{cb} \le \delta.
\end{equation}
Similarly, we have
\[
\|v_2x\| \le \tau(f(1-p)xx^*)^{1/2}
\]
hence by   \eqref{eq011} 
\begin{equation}\label{eq104}
\|v_2\|\le \|v_2\colon \ M\to R\|_{cb}\le \delta.
\end{equation}
We now turn to the following
\medskip

\n {\bf Claim 1.} $\vp>0$ can be chosen so that $\|\Lambda E(\vp)\|_2 \le 1$. 

For 
each integer $n\ge 1$, let
\[
\pi^n_2 = \sup\left\{\left(\sum\nolimits^n_1 \|\Lambda 
e_i\|^2\right)^{1/2}\right\}
\]
where the supremum runs over all possible orthonormal $n$-tuples $(e_1,\ldots, 
e_n)$ in the range of the projector $E(\vp)$. Note that for any operator
$w\colon
\ M\to E(\vp)$ and for  any $a_1,\ldots, a_n$ in  $M$ we have
\begin{equation}\label{eq105}
\left(\sum\nolimits^n_1 \|\Lambda w a_i\|^2\right)^{1/2} \le \pi^n_2\|w\| 
\left\|\sum 
a^*_ia_i\right\|^{1/2}.
\end{equation}
Indeed, let $T\colon \ [e_1,\ldots, e_n]\mapsto M$ be the map defined by $Te_i = 
a_i$. Note $\|T\|\le \|\Sigma a^*_ia_i\|^{1/2}$. 
 Let $F = \text{span}[wa_i]$
\begin{equation}
\sum \|\Lambda wa_i\|^2 = \sum^n_{i=1} \|\Lambda wTe_i\|^2 = \|\Lambda 
wT\|^2_2
\le \|\Lambda_{|F}\|^2_2 \|w\|^2 \|T\|^2
\end{equation}
and since $\dim F\le n$ we have $\|\Lambda_{|F}\|_2 \le \pi^n_2$ hence 
$\Sigma\|\Lambda wa_i\|^2 \le (\pi^n_2)^2 \|w\|^2 \|\Sigma a^*_ia_i\|$, which 
establishes \eqref{eq105}.
\medskip

\n {\bf Claim 2.} If $\pi^n_2 \le 1$, then we have
\begin{equation}\label{eq106}
\pi^n_2 \le (\vp\alpha   + \delta +\delta\pi^n_2)c.
\end{equation}
To prove this, consider $(e_1,\ldots, e_n)$ in $E(\vp)$ and let
\[
a = \sum^n_{i=1} e_{i1} \otimes (e_i, \Lambda e_i) \in M_n(G(\Lambda)).
\]
We have
\[
\|a\| = \max\left\{\left\|\sum e_{i1} \otimes e_i\right\|_{C_n\otimes_{\min} R}, 
\left\|\sum e_{i1} 
\otimes \Lambda e_i\right\|_{C_n\otimes_{\min} C} \right\}
\]
hence (since $\pi^n_2\le 1$)
\[
\|a\|  = \max\left\{1, \left(\sum \|\Lambda e_i\|^2\right)^{1/2}\right\} \le 1.
\]
By \eqref{eq101}, there is $\tilde a$ in $M_n(M)$ with $\|\tilde a\|\le c$, such 
that $(I\otimes u)(\tilde a) = a$. Clearly we may assume $\tilde a = \sum^n_1 
e_{i1} \otimes a_i$ with $a_i \in M$ such that $e_i = va_i$. Note that
\begin{equation}\label{eq107}
\left\|\sum a^*_ia_i\right\|^{1/2} = \|\tilde a\|\le c.
\end{equation}
Note that since $v = v_\delta +v_1+v_2$, we have $e_i = v_\delta a_i + v_1a_i + 
v_2a_i$, hence if we let $\Lambda_\vp = E(\vp)\Lambda = \Lambda E(\vp)$, we have
\begin{equation}\label{eq108}
\left(\sum\|\Lambda e_i\|^2\right)^{1/2} \le \left(\sum\|\Lambda_\vp v_\delta 
a_i\|^2\right)^{1/2} + 
\left(\sum \|\Lambda_\vp v_1a_i\|^2\right)^{1/2} + \left(\sum \|\Lambda_\vp 
v_2a_i\|^2\right)^{1/2}.
\end{equation}
By \eqref{eq102}, \eqref{eq107} and  \eqref{eq02} we have
\[
\left(\sum \|\Lambda_\vp v_\delta a_i\|^2\right)^{1/2} \le \vp \left(\sum
\|v_\delta  a_i\|^2\right)^{1/2} \le 
\vp(\alpha\tau(q))^{1/2}c
\]
and also by \eqref{eq103}
\[
\left(\sum \|\Lambda_\vp v_1a_i\|^2\right)^{1/2} \le \left(\sum \|\Lambda 
v_1a_i\|^2\right)^{1/2} \le 
\delta c.
\]
Finally, by \eqref{eq104}, \eqref{eq105} and \eqref{eq107}
 we have (recall $\Lambda_\vp =\Lambda E(\vp)$)
\[
\left(\sum\|\Lambda_\vp v_2a_i\|^2\right)^{1/2} \le \pi^n_2 \|E(\vp)v_2\| c\le 
\delta 
c\pi^n_2.
\]
Recapitulating, we can now deduce \eqref{eq106} from \eqref{eq108}, and Claim 2 
follows.

We can now prove Claim 1.

We assume $\vp < 1/2$. We will argue by contradiction. Assume that 
$\sup_m\pi^m_2>1$.  We will show that this is
impossible. Let 
$n+1$ be the smallest integer such that $\pi^{n+1}_2>1$. Note that $\pi^n_2\le 
1$ and $n\ge 1$ (because $\pi^1_2\le\vp <1$). Moreover, we have obviously 
$\pi^{n+1}_2 \le \pi^n_2 +\vp\le \pi^n_2 + 1/2$, hence $\pi^n_2 >1/2$. But now 
if 
we choose $\delta$ so that $\delta c < 1/2$, \eqref{eq106} implies
\[
\pi^n_2 \le c(\vp\alpha   +\delta) + (1/2)\pi^n_2
\]
hence
\[
\pi^n_2 \le 2c(\vp\alpha   + \delta),
\]
so that since $\pi^n_2>1/2$ we obtain
\[
1/2 \le 2c(\vp\alpha   + \delta).
\]
But now if we choose $\delta = 1/8c$ this implies
\begin{equation}\label{eq109}
1/4 \le 2c\vp\alpha  ,
\end{equation}
and here $\alpha=\alpha(\delta)$ is determined by $\delta$ but $\vp$ can still be made 
arbitrarily small. Thus we reach a contradiction, proving 
that $\sup\limits_m 
\pi^m_2\le 1$ for any $\vp<1/2$ for which \eqref{eq109} fails. This
proves  Claim 1 and completes the proof of Lemma 4.
\end{proof}

\begin{proof}[Proof of Theorem \ref{thm3}]
By Xu's result \eqref{eq100} we are reduced to $S$ of the form $S = G(\Lambda)$ 
for $\Lambda>0$ with dense range. By \eqref{eq1000}, we may assume that either 
$\Lambda$ or $\Lambda^{-1}$ has norm $\le 1$. But observe that if 
$\|\Lambda^{-1}\|\le 1$
\[
G(\Lambda^{-1}) = \{(x,\Lambda^{-1}x)\mid x\in C\} = \{(\Lambda y,y)\mid y\in 
\text{Dom}(\Lambda)\} \subset C\oplus R
\]
hence $G(\Lambda^{-1}) \simeq G(\Lambda)$ 
since the first is obtained from the second via the mapping
$(x,y) \to (y,x)$ which is obviously a complete  isometry
from $  C\oplus R  $ to   $  R\oplus C$.
 In 
particular, $G(\Lambda^{-1})$ embeds in $M_*$ iff 
$G(\Lambda)$ embeds in 
 $M_*$.
Thus to conclude we may as well assume that $\|\Lambda\|\le 1$. But then Lemma 
\ref{lem4} shows that for $\vp$ small enough we have a decomposition $R = H_\vp 
\oplus H^\bot_\vp$ and $\Lambda =\Lambda_\vp \oplus \Lambda'_\vp$ with 
$\|\Lambda_\vp\|_2 <\infty$. Clearly this implies $G(\Lambda) \simeq 
G(\Lambda_\vp) \oplus G(\Lambda'_\vp)$ but since $\|\Lambda_\vp\colon \ 
(H_\vp)_r \to C\|_{cb} = \|\Lambda_\vp\|_2<\infty$ we have $G(\Lambda_\vp) 
\simeq (H_\vp)_r$ and since $\vp \le \Lambda'_\vp \le 1$, we have obviously
(arguing as in the case when \eqref{eq(e)} holds) 
$G(\Lambda'_\vp )\simeq (H^\bot_\vp )_r \cap (H^\bot_\vp)_c$. This completes the 
proof of Theorem \ref{thm3}.
\end{proof}

\begin{rk} It may be worthwhile to point out that in Lemma \ref{lem4}, even if we know
that $G(\Lambda)^*$ is completely $c$-isomorphic to a subspace of $M_*$ with a fixed $c$, 
the $\vp$ given by Lemma \ref{lem4} may be arbitrarily small, and this happens
even for $M$ finite. Indeed, for the relevant examples, consider
a free circular    sequence  $(x'_i)$   on $(M,\tau)$ (with $\tau$ a normalized trace) and a
projection $p$ that is free from that family  and such that $\tau(p)=\vp$ (\cite{VDN}).
Then   $\overline{\text{span}}[p x'_i]$ provides the required phenomenon.
\end{rk}

In the next statement, we observe that Xu's decomposition for subspaces of 
$R\oplus C$ leads to an easy proof of  a result due to T.\ 
Oikhberg \cite{O} (with an improved bound), as follows.

\begin{thm}\label{thm5}
Let $S\subset R\oplus C$ be a closed subspace. If there is a completely bounded 
projection $P\colon \ R\oplus C\to S$ then there are Hilbert spaces $H,K$ such 
that $S\simeq H_r\oplus K_c$. Moreover there is a numerical constant $C$ such 
that $d_{cb}(S, H_r\oplus K_c) \le C\|P\|_{cb}$.
\end{thm}

\begin{proof}
By Xu's decomposition and the above remarks, it suffices to prove this for $S = 
G(\Lambda)$ with $0<\Lambda$ and $\|\Lambda\|\le 1$. Then the projection $P$ can 
be written as
\[
\forall (x,y) \in R\oplus C\qquad P(x,y) = (\alpha x + \beta y, \Lambda(\alpha x 
+ \beta  y))
\]
where $\alpha\in CB(R,R)$ and $\beta\in CB(C,R)$. By restricting $P$, we find
\begin{equation}\label{eq110}
\max\{\|\alpha\|_{CB(R,R)}, \|\Lambda\alpha\|_{CB(R,C)} \|\beta\|_{CB(C,R)}, 
\|\Lambda\beta\|_{CB(C,C)}\}\le \|P\|_{cb}.
\end{equation}
Moreover since $P$ is a projection onto $G(\Lambda)$ we have for any $x$ in $R$
\[
\alpha x + \beta\Lambda x = x
\]
hence
\[
\Lambda \alpha + \Lambda \beta \Lambda = \Lambda
\]
which implies by \eqref{RC2} and   \eqref{eq110} (since we assume $\Lambda \le 1$)
\[
\|\Lambda\|_{CB(R,C)} =\|\Lambda\|_{2}\le
\|\Lambda\alpha\|_{2} + 
\|\Lambda\beta\Lambda\|_{2}=
 \|\Lambda\alpha\|_{CB(R,C)} + 
\|\Lambda\beta\Lambda\|_{CB(C,R)} \le 2\|P\|_{cb}.
\]
Thus we conclude
\[
\|\Lambda\|_{CB (R,C)} \le 2\|P\|_{cb}
\]
and hence the map $u\colon \ x\to (x,\Lambda x)$ is a complete isomorphism 
between $R$ and $G(\Lambda)$ with
$$
d_{cb}(R,G(\Lambda)) \le \|u\|_{cb} \|u^{-1}\|_{cb} \le 2\|P\|_{cb}.\qquad \qed
$$
\renewcommand{\qed}{}\end{proof}

\begin{rk}

The preceding statement yields a rather satisfactory estimate in the following 
result from \cite{PS}:\ If an operator space $E$ is exact as well as its dual, 
then there are Hilbert spaces, $H,K$ such that $E\simeq H_r\oplus K_c$ and 
moreover 
\[
d_{cb}(E,H_r\oplus K_c) \le 2^{5/2} ex(E) ex(E^*)
\]
where $ex(E)$ denotes the exactness constant of $E$. This seems rather sharp 
when $ex(E) ex(E^*)$ is large:\ Consider for instance the case $E = OH_n$, we 
have then (cf.\ \cite[p. 336]{P0}) $ex(E) = ex(E^*) \simeq n^{1/4}$ but on the 
other hand it is easily checked that 
\[
d_{cb}(OH_n, H_r\oplus K_c) \simeq n^{1/2} \simeq ex(E) ex(E^*).
\]

\end{rk}

\bigskip\bigskip

\setcounter{thm}{0}

\centerline{\bf Appendix}

\bigskip
In this appendix, we will reprove  Junge's result \cite{J} that $OH$
embeds completely  isomorphically   into a non-commutative
$L_1$-space. The main  idea is the same as his, but our exposition is
shorter and makes more  transparent the relationship between the methods
from \cite{J} and
\cite{PS}.  We base the argument on the complex interpolation method
instead of the  Pusz--Woronowicz formula. Actually, there is nothing
mysterious there:\ indeed  the ``purification of states'' associated in
\cite{PW1} (see also \cite{PW2,W1,W2,W3})  to a pair of faithful  states
$(\varphi,\psi)$ on a
$C^*$-algebra $A$ is known to be very closely  related to the complex
interpolation space $(A_0,A_1)_{1/2}$ where the Hilbert  spaces $A_0,A_1$
are obtained by completing $A$ for the norms
\[
\|x\|_{A_0}  =  (\varphi(x^*x))^{1/2},\qquad \|x\|_{A_1} = (\psi(xx^*))^{1/2}.
\]
This close connection has been explored
 in depth notably by B.\ Simon, Uhlman,
Peetre and probably others, besides Pusz and Woronowicz.

The proof rests on the following basic fact which had been known
to the author (and probably also to Junge)  for some
time, before Junge proved his
embedding result for $OH$. A detailed proof is included as the solution to
Exercise 7.8 in \cite{P3}.
We  reproduce it here for the convenience of the reader.

\begin{pp}\label{proA1}
$OH$ is completely isometric to a quotient of a subspace of $R\oplus C$.
\end{pp}

\begin{proof}
 Let $\mu$ be the harmonic measure of
the  point $z = 1/2$ in the strip $S = \{z\in \CC\mid 0 < \hbox{Re}(z) <
1\}$. Recall  that $\mu$ is a probability measure on $\partial S$ such
that $f(1/2) = \int f\  d\mu$ whenever $f$ is a bounded harmonic function
on $S$ extended  non-tangentially to $\overline S$. Obviously $\mu$ can
be written as $\mu =  2^{-1}(\mu_0 + \mu_1)$ where $\mu_0$ and $\mu_1$
are probability measures  supported respectively by
$$\partial_0 = \{z\mid \hbox{Re}(z) = 0\}\quad \hbox{and}\quad \partial_1 =
\{z\mid \hbox{Re}(z) = 1\}.$$
Let $(A_0,A_1)$ be a compatible pair of Banach spaces. We first need to
describe
$(A_0,A_1)_{1/2}$ as a quotient of a subspace of $L_2(\mu_0;A_0)\oplus
L_2(\mu_1;A_1)$. The classical argument for this is as follows.

 We denote by ${\cal F}(E_0,E_1)$ the set of all
bounded continuous functions $f\colon \ \ovl S\to E_0+E_1$ which are
holomorphic on $S$ and such that $f_{|\partial_0}$ and $f_{|\partial_1}$
are bounded continuous functions with values respectively in $E_0$ and
$E_1$.

  We start by
showing that for any $x$ in $(A_0,A_1)_{1/2}$ we have
$$\|x\|_{(A_0,A_1)_{1/2}} = \inf\{\max\{\|f\|_{L_2(\mu_0;A_0)},
\|f\|_{L_2(\mu_1;A_1)}\}$$
where the infimum runs over all $f$ in ${\cal F}(A_0,A_1)$ such that
$f(1/2) = x$.
 For a
proof, see      \cite[p.~224]{KPS}. Let  then
 $E = L_2(\mu_0; A_0)\oplus_\infty
L_2(\mu_1,A_1)$ and let $G\subset E$ be the closure of the subspace
$\{f_{|\partial_0} \oplus f_{|\partial_1}\mid f\in {\cal F}(A_0,A_1)\}$.
The preceding  equality shows that the mapping $f\to f(1/2)$ defines a
metric surjection
$Q\colon \ G\to (A_0,A_1)_{1/2}$. We now consider the couple $(A_0,A_1) =
(R,C)$, where we think of $R$ and $C$ as operator space stuctures
on the ``same"
  underlying  vector space, identified with $\ell_2$. We introduce the
operator space
$E = L_2(\mu_0;\ell_2)_r
\oplus  L_2(\mu_1;\ell_2)_c$. Let $G$ and $Q\colon \ G\to \ell_2$ be the
same as before. Note that $G$ is nothing but 
the $\ell_2$-valued version of the Hardy space $H^2$ 
on the strip $S$,  
 so that if we assume $f$ analytically extended inside
$S$, we have $Q(f) =  f(1/2)$.

\n We first claim that
$$\|Q\colon \ G\to OH\|_{cb} \le 1.$$
To verify this, consider $x$ in $M_n(G)$ with $\|x\|_{M_n(G)} \le 1$. We claim
that $\|x(1/2)\|_{M_n(OH)}\le 1$. We may view $x$ as a sequence $(x_k)$ of
$M_n$-valued functions on $\partial S$ extended analytically inside $S$, so
that
$$\|x\|_{M_n(G)} = \max\left\{\left\|\left(\int \sum x_kx^*_k\
d\mu_0\right)^{1/2}
\right\|_{M_n},\quad \left\|\left(\int \sum x^*_kx_k\
d\mu_1\right)^{1/2}\right\|_{M_n}\right\},$$
and by   \cite[(7.3)$'$]{P3}
$$\|x(1/2)\|_{M_n(OH)}^2 = \left\|\sum x_k(1/2) \otimes
\ovl{x_k(1/2)}\right\|_{\rm  min}  = \sup\left\{\left|\hbox{tr}\left(\sum
x_k(1/2) ax_k(1/2)^* b\right)\right|\right\}$$
where the supremum runs over all $a,b\ge 0$ in $M_n$ such that $\hbox{tr}|a|^2
\le 1$ and $\hbox{tr}|b|^2 \le 1$. Fix $a,b$ satisfying these conditions.
Consider then the analytic function
$$F(z) = \hbox{tr}\left(\sum x_k(z) a^{2z}x_k(\bar z)^*
b^{2(1-z)}\right),$$ on $S$. Note that
$$F(1/2) = \hbox{tr}\left( \sum x_k(1/2) a x_k(1/2)^* b\right) =
2^{-1}\left(
\int_{\partial_0} F\ d\mu_0 + \int_{\partial_1} F\ d\mu_1\right).$$
But for all $z = it$ in $\partial_0$ we have
$$F(it) = \sum_k \hbox{tr}(b^{1-it} x_k(it)a^{2it} x_k(-it)^* b^{1-it})$$
hence by Cauchy--Schwarz for any $z$ in $\partial_0$
$$|F(z)| \le \left(\sum_k \hbox{tr}(bx_k(z) x_k(z)^*b)\right)^{1/2}
\left(\sum_k
\hbox{tr}(b x_k(\bar z)x_k(\bar z)^* b)\right)^{1/2}.$$
A similar verification shows that for any $z$ in $\partial_1$ we have
$$|F(z)| \le \left(\sum_k \hbox{ tr}(ax_k(z)^* x_k(z)a)\right)^{1/2}
\left(\sum_k \hbox{ tr}(ax_k(\bar z)^* x_k(\bar z)a)\right)^{1/2}.$$
Thus we obtain by Cauchy--Schwarz
$$ {|F(1/2)|=|\int F d\mu| \le 2^{-1}\left( \int_{\partial_0} |F|
\ d\mu_0 +
\int_{\partial_1} |F|\  d\mu_1\right)}  $$
$$\le 2^{-1}\left\{
\hbox{tr}\left(b^2 \int \sum x_k x^*_k \ d\mu_0\right)
\quad + \hbox{tr}\left(a^2 \int \sum x^*_kx_k \ d\mu_1\right)
\right\}\le
\|x\|_{M_n(G)}\le 1,$$ which proves our claim.

It is now easy to show that $Q$ is actually a complete metric surjection, or
equivalently, that $I\otimes Q \colon \ M_n(G) \to M_n(OH)$ is a metric
surjection for any $n\ge 1$. Indeed, consider $x\in M_n(OH)$ with
$\|x\|_{M_n(OH)} < 1$. Since  $M_n(OH) = (M_n(R),
M_n(C))_{1/2}$ (isometrically) by   \cite[Corollary 5.9]{P3},  there is a
bounded continuous analytic function $f$ on
$\ovl S$ with values in
$M_n(R) + M_n(C)$ such that $$\alpha_0=\sup\{\|f(z)\|_{M_n(R)}\mid z\in
\partial_0\} < 1, \
 \alpha_1=\sup\{\|f(z)\|_{M_n(C)} \mid z\in \partial_1\} < 1\ {\rm
and}\  f(1/2) = x.$$ Let us  write $f(z) = (f_k(z))_k$ where $f_k$ is
an
$M_n$-valued function on $\ovl S$.  We have trivially
$$\left\|\left(\int \sum f_k(z) f_k(z)^* d\mu_0(z)\right)^{1/2}\right\|_{M_n}
\le \alpha_0 < 1$$
and
$$\left\|\left(\int \sum f_k(z)^* f_k(z) d\mu_1(z)\right)^{1/2}\right\|_{M_n}
\le \alpha_1 < 1$$
hence $\|f\|_{M_n(G)} < 1$. Since clearly $(I\otimes Q)(f) = x$, this shows
that
$I\otimes Q\colon \ M_n(G)\to M_n(OH)$ is a metric surjection. Thus we have
completely isometrically $OH\simeq G/\ker(Q)$.
Finally since $G\subset R\oplus C$ this completes the proof.
\end{proof}

Let $E = L_2(\mu_0;\ell_2)_r\oplus L_2(\mu_1;\ell_2)_c$. The preceding
argument
shows that
\begin{equation}\label{G}
OH \simeq G/N
\end{equation}
where $G\subset E$ is the subspace of boundary values of analytic functions on
the strip $S = \{0< \text{Re } z < 1\}$, and where $N$ is the subspace of $G$
formed of all $f$ in $G$ such that $f(1/2) = 0$. Thus, $OH$ appears as a
quotient, namely $G/N$, of a subspace, namely $G$, of $R\oplus C$ since
obviously
$E \simeq R\oplus C$. Moreover, the subspace $G\subset E$ is the {\em graph\/}
of a (necessarily closed) unbounded operator $T\colon \ \text{Dom}(T) \to
L_2(\mu_1; \ell_2)_c$ where $\text{Dom}(T) \subset L_2(\mu_0,\ell_2)_r$ is the
dense subspace formed of all the restrictions $f_{|\partial_0}$ when $f$ runs
over $G$. Since $G$ is formed of {\em analytic\/} functions, the
restriction of
$f$ to  $\partial_0$ (or $\partial_1$) entirely determines $f$, therefore
$f\in
G\to f_{|\partial_0}$ and $f\in G\to f_{|\partial_1}$ are one to one, so that
the definition of $T$ is clear:\ we simply have
\[
T(f_{|\partial_0}) = f_{|\partial_1}.
\]
Note that $T$ has dense range. By the polar decomposition of $T$ (cf.\
\cite[p.~1249]{DS}) we have $T = U|T|$ where $U\colon \ L_2(\mu_0;\ell_2)\to
L_2(\mu_1;\ell_2)$ is unitary and where $| T|\colon \ L_2(\mu_0;\ell_2) \to
L_2(\mu_1;\ell_2)$ is an unbounded, $\ge 0$ and self-adjoint operator.

Clearly, since $L_2(\mu_0;\ell_2)_r\simeq R$ and $L_2(\mu_1; \ell_2)_c
\simeq C$
are ``homogeneous'' operator spaces (i.e.\ for any $u\colon \ R\to R$ or
$u\colon \ C\to C$ we have $\|u\|_{cb} = \|u\|$), $U$ 
(or its inverse) is completely isometric
from $L_2(\mu_0;\ell_2)_c$ to $L_2(\mu_1;\ell_2)_c$, and hence
$  I\oplus U^{-1} $ is  completely isometric on $L_2(\mu_0,\ell_2)_r \oplus L_2(\mu_1,\ell_2)_c$.
 Let $\Lambda\colon \
L_2(\mu_1;\ell_2)_r\to L_2(\mu_1;\ell_2)_c$ be the same map as $|T|$ but
viewed
as acting between the indicated operator spaces (so that $T=U\Lambda$). 

Then we have trivially
\[
G\simeq (I\oplus U^{-1})(G) 
\]
but
\[
(I\oplus U^{-1})(G)= \{(x, \Lambda x) \mid x \in
\text{Dom}(\Lambda)\}.
\]
So we are reduced to the following.

\begin{pp}\label{proA2}
Let $\Lambda\colon \ R\to C$ be a closed self-adjoint densely defined
unbounded
operator with $\Lambda\ge 0$. Let
\[
G(\Lambda) = \{(x,\Lambda x)\mid x\in \text{Dom}(\Lambda)\} \subset R\oplus C
\]
be the graph of $\Lambda$. Then the dual $G(\Lambda)^*$
 is completely 2-isomorphic to a subspace of   a non-commutative
$L_1$-space. In particular $OH$ embeds completely 2-isomorphically
in $M_*$ for some von Neumann algebra $M$.
\end{pp}

\begin{proof}
Let $\{E_\alpha\}$ be a net of finite dimensional subspaces of 
$\text{Dom}(\Lambda)$ directed by inclusion and such that $\cup E_\alpha = 
\text{Dom}(\Lambda)$. Let $G_\alpha = \{(x,\Lambda x)\mid x\in E_\alpha\}$. Then 
$G(\Lambda) = \cup G_\alpha$ (directed union) and hence for any c.b.\ map 
$u\colon \ G(\Lambda)\to M_n$ we have
\[
\|u\|_{cb} = \lim \uparrow\|u_{|G_\alpha}\colon \ G_\alpha \to M_n\|_{cb}.
\]
It follows that $G(\Lambda)^*$ embeds completely isometrically into an 
ultraproduct of the spaces $G^*_\alpha$. Since by \cite{R}, ultraproducts 
preserve the class of subspaces of non-commutative $L_1$-spaces
(the operator space version of this is easy to derive from
\cite{R}) we are reduced 
  to proving  this with
$G(\Lambda)$  replaced by $G_\alpha$. In that case we may as well replace
$C$ by
$C_n$ (where $=\dim(E_\alpha)$) and replace $R$ by $R_n$.

\n  Thus we are reduced to proving
the result for $G(\Lambda)\subset R_n\oplus C_n$ for some invertible operator
$\Lambda\ge 0$ from $R_n$ to $C_n$.
In that case, we may as well assume (by homogeneity) that $\Lambda e_{1i}
=
\lambda_i e_{i1}$ for  some $\lambda_i>0$. But then this  follows
from the next result which is  somewhat implicit in Shlyakhtenko's work
\cite{S1}, and in any case is included in the above Lemma \ref{lem22}.
\end{proof}

\begin{pp}\label{proA3}
With the notation as in the first part. Let $I = \{1,\ldots, n\}$ and
\[
a_i = \ell_i + \lambda_i\ell'{}^*_i.
\]
Let $G_n = \text{span}[a_1,\ldots, a_n]$ and let $W_n$ be the von~Neumann
algebra generated by $G_n$. Then $G_n$ is  completely 2-isomorphic to
$G(\Lambda_n)\subset R_n\oplus C_n$ and $G_n$ is completely 2-complemented in
 $W_n$. More precisely, we have a surjective mapping
$P_n:\ W_n \to G(\Lambda_n)$ with $\|P_n\|_{cb}\le 1$ such that
$G(\Lambda_n)$ is   completely 2-isomorphic to the quotient
$W_n/\ker(P_n)$.
Therefore, $G(\Lambda_n)^*$
is completely 2-isomorphic to a subspace of 
a non-commutative $L_1$-space, namely   the predual of  $W_n$.
 \end{pp}

\begin{proof}
We let $P$ be as in the proof of Lemma \ref{lem22}. Let 
$V:\ G_n\to G(\Lambda_n)$ be defined by
$V(\ell_i + \lambda_i\ell'{}^*_i)= e_{1i} \oplus \lambda_i e_{i1}$.
Finally, let $P_n=VP:\ W_n\to G(\Lambda_n)$. The proof of Lemma \ref{lem22} shows that
$\|P_n\|_{cb}\le 1$ and by the triangle inequality  we have $\|V^{-1}\|_{cb}\le 2$. 
Therefore, $  G(\Lambda_n)$ is completely 2-isomorphic to $W_n/\ker(P_n)$.
\end{proof}
\noindent{\bf Note:} In the above we used a discretization of $\Lambda$
to make the proof as elementary as possible, but this is not really
necessary if one uses the general picture described in \cite{S1}.
This alternate route is much more elegant but perhaps a bit more ``abstract''. We will 
merely outline it. We consider the (complex) Hilbert space $H= 
L_2(\mu_0;\ell_2)\oplus L_2(\mu_1,\ell_2)$ equipped with the norm
\[
\|(x_0,x_1)\| = (2^{-1}(\|x_0\|^2 + \|x_1\|^2))^{1/2}.
\]
As is classical, any $x = (x_0,x_1)$ admits (via Poisson integrals) a harmonic 
extension inside $S$, i.e.\ there is a harmonic function $\tilde x\colon \ S\to 
\ell_2$ such that $\|\tilde x(\cdot)\|^2$ admits a harmonic majorant $u$ on $S$ 
and admitting $x_0$ and $x_1$ as its nontangential boundary values respectively 
on $\partial_0$ and on $\partial_1$. Note that $\|x\| = \inf\{u(1/2)\}$ where 
the infimum runs over all such majorants and the Poisson integral of the 
function equal to $\|x_0(\cdot)\|^2$ on $\partial_0$ and $\|x_1(\cdot)\|^2$ on 
$\partial_1$ produces the minimal $u$.

We will denote by $h^2(\ell_2)$ the space of all harmonic functions $\tilde x$ 
obtained in this way. All such functions are implicitly extended
 nontangentially to the closure of $S$.
 Thus $h^2(\ell_2)$ can be identified with $H$. We denote 
by $H^2(\ell_2)$ the subspace of $h^2(\ell_2)$ formed of all the {\em 
analytic\/} functions. The spaces $h^2(\ell_2)$ and $H^2(\ell_2)$ may be viewed
as conformally equivalent copies of the usual spaces on the unit disc.

 \n For any $f=(f_k)$ in $h^2(\ell_2)$, we set $\bar f=(\bar f_k)$.

\n We denote by $H_{\bb R}$ the real linear subspace of $H$ 
of all elements of the form $(\bar f_{|\partial_0},  f_{|\partial_1})$ when $f$ 
runs over all functions in $H^2(\ell_2)$. Note that the map
\[
j\colon \ H^2(\ell_2)\to (\bar f_{|\partial_0}, f_{|\partial_1})\in H
\]
is a {\it real linear} isometry, with range $H_{\bb R}$. It is easy to
check that
$H_{\bb R}\cap i H_{\bb R} = 
\{0\}$ (because an analytic function in $H^2(\ell_2)$ that vanishes
on $\partial_1$ must vanish everywhere) 
and that $H_{\bb R} +iH_{\bb R}$ is dense in $H$
(because if an element of $h^2(\ell_2)$ is supported on $\partial_1$ or 
$\partial_0$ and is orthogonal to $H^2(\ell_2)$, it must be anti-analytic, and hence
must vanish identically; therefore
 the restrictions $\{ f_{|\partial_1} \mid f\in H^2(\ell_2)\} $
are dense in $L_2(\mu_1;\ell_2)$, and similarly for $\partial_0$).

As pointed out in \cite[Remark 2.6]{S1}, the basic construction of
\cite{S1} can be carried out starting from the data
of the embedding $j  :\ H_{\bb R} \to H$, using \cite {RV} to obtain 
a  group of orthogonal transformations of $H$
satisfying the KMS condition relative to this embedding.
Let $\cl F$ be the full Fock space over $H$.
We will identify $H$ with $L_2(\partial_0 \cup \partial_1; \mu)$.
 With the previous notation
we set for any $f$ in $H^2(\ell_2)$
\begin{equation}\label{t}
  t(f) =(\ell(  \bar f 1_{\partial_0}   ))^* + \ell(    f 1_{\partial_1}   ).
\end{equation}
Observe that $f\to t(f)$ is now a {\it complex linear} isomorphic embedding of $H^2(\ell_2)$ into $B(\cl F   )$.
Note that this ``quantization" of $H^2(\ell_2)$ 
seems to  be of independent interest (even 
for scalar valued Hardy spaces, when $\ell_2$ is replaced by
${\bb C}$). More generally, \eqref{t} makes senses for any $f$ in
$h^2(\ell_2)$; the resulting mapping is then a completely isomorphic
embedding of $L_2(\mu_0;\ell_2)_r\oplus L_2(\mu_1;\ell_2)_c$ into $B(\cl
F   )$.

Shlyakhtenko \cite{S1} made an extensive study of the 
von Neumann algebra $M$ generated by the operators
$\{s(h)=\ell(h) + \ell(h)^*,\ h\in H_{\bb R}\}$. Since for any $f$ in $H^2(\ell_2)$,
$$2t(f) =s(j(f)) -i s(j(if)),$$ $M$ is  equivalently generated
by the family $\{t(f), \ f\in H^2(\ell_2)\}$.
Finally, arguing as for the above Lemma \ref{lem22}, we see that there
is a projection $P:\ M \to t(H^2(\ell_2))$ with $\|P\|_{cb}\le 2$, and
$t(H^2(\ell_2))$ is completely isomorphic to the
space $G$ appearing in \eqref{G}. Thus we can conclude, as  
in  Proposition A.\ref{proA3}, that
$G$ (and a fortiori $OH$) is completely 2-isomorphic to a quotient of
$M$, via a normal surjection $M\to OH$.  Thus, taking adjoints,
we find that $OH$ embeds completely isomorphically  into $M_*$.
\begin{rk}
The constant $2$ appearing in Proposition A\ref{proA2} is better than in
Junge's
original work. This is not too significant but the question whether
$OH$ embeds completely isometrically into $M_*$ for some von Neumann
algebra $M$ remains open. 
\end{rk}
\begin{rk} Our proof does not yield  the fact  announced by Junge
(yet unpublished) that, in the above Proposition  A.\ref{proA2}, $M$ can be chosen hyperfinite.
 Note however, that since by \cite{PS} we   obtain an $M$ that is 
a quotient of a $C^*$-algebra with the $WEP$ ($QWEP$), one can deduce
from the strong local reflexivity of non-commutative $L_1$-spaces 
(see \cite{EJR}) Junge's result that
for each $n$ and $c>2$ there is an integer $N$ and a subspace $E_n\subset S_1^N$ such that
$d_{cb}(E_n,OH_n)\le c$.
\end{rk}

\begin{rk}
The same proof suitably modified shows that $OH$ embeds
 completely isomorphically 
 in a  non-commutative $L_p$-space for any $p$ with $1<p<2$. (The case
$p=2$ is of  course trivial.) That result was known to Junge and Xu. 
Indeed, for any $0<\theta<1$, we have by \cite{P1}
$OH = (R[\theta], R[1-\theta])_{\frac12}$, hence (arguing as for Proposition 
A.\ref{proA1}) we find that $OH$ is a quotient of a subspace of $R[\theta]\oplus 
R[1-\theta]$.  Now let $p = (1-\theta)^{-1}$ as before. In that case we claim 
that $(R[\theta]\oplus R[1-\theta])^*$ embeds completely isomorphically 
into $S_p$. Indeed, as we  already mentioned, $R[\theta]^* = R[1-\theta]$
(resp.\ $R[1-\theta]^* =  R[\theta]$) can be identified with the subspace
of column (resp.\ row) matrices  in $S_p$. This proves our claim.

\n  More generally, it follows from Xu's results in \cite{X,X2} 
(see also \cite{JX} for related facts) that for any closed 
unbounded positive operator $\Lambda\colon \ R[\theta]\to R[1-\theta]$ with 
dense domain, dense range and zero kernel, the graph $G(\Lambda) \subset 
R[\theta] \oplus R[1-\theta]$ is such that $G(\Lambda)^*$ embeds in a 
non-commutative $L_p$-space. Thus by the same principle as above for $p=1$, we 
can show that $OH$ embeds completely isomorphically  in a non-commutative
$L_p$ for all $1<p<2$. See \cite{X2} for more on this theme.
\end{rk}

\begin{rk}
Junge observed already in \cite{JO} that $OH$ does not embed completely isomorphically 
into a    non-commutative $L_q$-space for $2<q<\infty$. Actually, in that
case it is even  impossible to embed $OH_n$ uniformly over $n$ into such
a space. For the reader's convenience, we now sketch Junge's argument for
this fact. We will use the 
non-commutative Khintchine inequalities due to Lust--Piquard (cf.\ 
\cite[p.~193]{P3}). For our present purpose, it is convenient to state them 
using 
the ``vector valued'' version of the Schatten classes introduced in \cite{P2} 
and denoted by $S_p[E]$ where $1\le p<\infty$ and $E$ is an arbitrary operator 
space.
Let $(\vp_k)$ denote the Rademacher functions on $(\Omega,m)$ where (say) 
$\Omega= [0,1]$ and $m$ is normalized Lebesgue measure. Then if the operator 
space $E$ is assumed to be (completely isometrically) a subspace of a 
non-commutative $L_q$-space and if $2\le q<\infty$, then for any finite sequence 
$a_1,\ldots, a_n$ in $E$ we have
\begin{equation}\label{eq(k)}
\left(\int\left\|\sum\nolimits^n_1 \vp_ka_k\right\|^q dm\right)^{1/q} \le B_q 
\left(\left\|\sum\nolimits^n_1 e_{1k} \otimes a_k\right\|_{S_q[E]} + 
\left\|\sum\nolimits^n_1 e_{k1}\otimes a_k\right\|_{S_q[E]}\right)
\end{equation}
where $B_q$ is a constant depending only on $q$. See \cite{JX} for the extension 
of \eqref{eq(k)} to the case of general non-commutative $L_q$-spaces, including 
the non-semifinite case. 

\n Now, let $u\colon \ OH_n\to E$ be a linear isomorphism 
and let $(e_1,\ldots, e_n)$ denote an orthonormal basis in $OH_n$. We have 
clearly
\[
\sqrt n \le \|u^{-1}\| \left(\int\left\|\sum\nolimits^n_1 \vp_kue_k\right\|^q 
dm\right)^{1/q}.
\]
On the other hand  applying \eqref{eq(k)} to $a_k=ue_k$ and using 
\cite[Cor.~1.2]{P2} we find 
\[
\left(\int\left\|\sum\nolimits^n \vp_kue_k\right\|^q dm\right)^{1/q} \le 
\|u\|_{cb} B_q \left(\left\|\sum\nolimits^n_1 e_{1k} \otimes 
e_k\right\|_{S_q[OH_n]} + \left\|\sum\nolimits^n_1 e_{k1} \otimes 
e_k\right\|_{S_q[OH_n]}\right).
\]
Finally, by easy interpolation arguments (based on \cite[Cor.~1.4]{P2}) we find
\[
\left\|\sum\nolimits^k_1 e_{1k}\otimes e_k\right\|_{S_q[OH_n]} \le n^{\frac1{2q} 
+ \frac14}
\]
and similarly for $\left\|\sum^k_1 e_{k1} \otimes e_k\right\|_{S_q[OH_n]}$. Thus 
we conclude
\[
n^{\frac12} \le \|u^{-1}\|\|u\|_{cb} 2B_q n^{\frac1{2q}+\frac14}
\]
and a fortiori we find
\[
d_{cb}(E,OH_n) = \inf\|u\|_{cb} \|u^{-1}\|_{cb} \ge (2B_q)^{-1} 
n^{\frac14-\frac1{2q}}.
\]
A similar argument can be applied with $(C_n,R_n)_\theta$ instead of $OH_n$. The 
same calculations yield that for any $q> \max\{p,p'\}$ with $p = 
(1-\theta)^{-1}$ and $p' =\theta^{-1}$, we have
\[
d_{cb} (E, (C_n,R_n)_\theta) \ge (2B_q)^{-1} n^{\beta/2}
\]
where $\beta = \min\{p^{-1}-q^{-1}, p'{}^{-1}-q^{-1}\}$.
Note however that if $p = (1-\theta)^{-1}$ then $(C_n,R_n)_\theta$ embeds 
completely isometrically in both $S_p$ and $S_{p'}$; indeed it can be identified 
with $\text{span}[e_{1k}\mid 1\le k\le n]$ in $S_p$ and with 
$\text{span}[e_{k1}\mid 1\le k\le n]$ in $S_{p'}$.
\end{rk}

\noindent {\bf Acknowledgement:}\ I am grateful to Quanhua Xu for fruitful
conversations, useful correspondence and for reminding 
me of some results that I should not have
forgotten. Thanks also to the referee for remarks
that led to an improved   presentation.

\vfill\eject

\end{document}